\documentclass[a4paper]{amsart}

\usepackage{amssymb,amsmath,graphics}
\usepackage{latexsym}
\usepackage{eucal}
\makeatletter
\@addtoreset{equation}{section}\makeatother

\newtheorem{theo}{Theorem}[section]
\newtheorem{defi}[theo]{Definition}
\newtheorem{lem}[theo]{Lemma}
\newtheorem{prop}[theo]{Proposition}
\newtheorem{cor}[theo]{Corollary}
\theoremstyle{remark} \newtheorem{remark}[theo]{Remark}
\newcommand{\mc}{\mathcal}
\newcommand{\rr}{\mathbb{R}}
\newcommand{\nn}{\mathbb{N}}
\newcommand{\cc}{\mathbb{C}}
\newcommand{\hh}{\mathbb{H}}
\newcommand{\zz}{\mathbb{Z}}

\newcommand{\la}{\lambda}
\newcommand{\eps}{\epsilon}

\newcommand{\psdo}{\Psi\textrm{DO}}
\newcommand{\pl}{\partial}
\newcommand{\x}{\times}

\newcommand{\til}{\widetilde}

\newcommand{\supp}{\textrm{supp}}
\newcommand{\cjd}{\rangle}
\newcommand{\cjg}{\langle}

\newcommand{\demi}{\frac{1}{2}}
\newcommand{\ndemi}{\frac{n}{2}}
\newcommand{\tra}{\textrm{Tr}}

\newcommand{\diag}{\textrm{diag}}

\newcommand{\rang}{\mathop{\hbox{\rm rank}}\nolimits}

\newcommand{\trans}{{^t}\!}
\newcommand{\TR}{\mathop{\hbox{\rm TR}}\nolimits}
\newcommand{\odd}{\textrm{odd}}
\newcommand{\sing}{\textrm{sing}}
\newcommand{\reg}{\textrm{reg}}
\newcommand{\indic}{\operatorname{1\negthinspace l}}
\newcommand{\zerotr}{\mathop{\hbox{\rm 0-Tr}}\nolimits}
\newcommand{\zerov}{\mathop{\hbox{\rm 0-vol}}\nolimits}

\def\qed{\hfill$\square$}
\begin{document}
\title[Generalized Krein formula, determinants and Selberg zeta function]{Generalized Krein formula, determinants and Selberg zeta function in even dimension}
\author{Colin Guillarmou}
\address{Laboratoire J.-A. Dieudonn\'e\\
U.M.R. 6621 du C.N.R.S.\\
Universit\'e de Nice\\
Parc Valrose\\
06108 Nice Cedex 02\\
France}
     \email{cguillar@math.unice.fr}

\subjclass[2000]{Primary 58J50, Secondary 47A40, 11M36, 58J52}
%
\begin{abstract}
\noindent For a class of even dimensional asymptotically hyperbolic (AH) manifolds,
we develop a generalized Birman-Krein theory to study scattering asymptotics and, when the curvature is constant, 
to analyze Selberg zeta function. 
The main objects we construct for an AH manifold $(X,g)$ are, on the first hand, 
a natural spectral function $\xi$ for the Laplacian $\Delta_g$, which replaces the counting
function of the eigenvalues in this infinite volume case, and on the other hand
the determinant of the scattering operator $S_X(\la)$ of $\Delta_g$ on $X$.
Both need to be defined through regularized functional: renormalized trace on the bulk $X$ 
and regularized determinant on the conformal infinity $(\pl\bar{X},[h_0])$.
We show that $\det S_X(\la)$ is meromorphic in $\la$, with divisors given 
by resonance multiplicities and dimensions of kernels of GJMS conformal Laplacians 
$(P_k)_{k\in\nn}$ of $(\pl\bar{X},[h_0])$, moreover $\xi(z)$ is proved to be the phase of $\det S_X(\ndemi+iz)$ 
on the essential spectrum $\{z\in\rr^+\}$. 
Applying this theory to convex co-compact quotients $X=\Gamma\backslash\hh^{n+1}$ of hyperbolic space $\hh^{n+1}$, 
we obtain the functional equation $Z(\la)/Z(n-\la)=(\det S_{\hh^{n+1}}(\la))^{\chi(X)}/\det S_X(\la)$ 
for Selberg zeta function $Z(\la)$ of $X$, where $\chi(X)$ is the Euler characteristic of $X$. This 
describes the poles and zeros of $Z(\la)$, computes $\det P_k$ in term of $Z(\ndemi-k)/Z(\ndemi+k)$ 
and implies a sharp Weyl asymptotic for $\xi(z)$.
\end{abstract}
\maketitle

\section{Introduction}

The study of Selberg's zeta function $Z(\la)$ for hyperbolic manifolds has lead to many 
fascinating results relating dynamic and spectral-scattering theory of the Laplacian. 
The function $Z(\la)$ (for $\la\in\cc,\Re(\la)\gg 0$) is defined as a convergent infinite product 
over the closed geodesics of the hyperbolic manifold $X:=\Gamma\backslash\hh^{n+1}$, 
in other words it is a purely dynamical function.   
For compact or finite volume hyperbolic manifolds, it has been proved 
that it extends meromorphically to $\cc$ and, more interestingly, its zeros and poles 
are given by topological data on the one hand and by spectral/scattering data of the Laplacian $\Delta_X$
on the other hand (see for instance \cite{Se,GA,GW,He,He2,LP,BO1}).
This may be considered as a direct relation between classical dynamics and its quantization. 

It is worth to emphasize that, unlike for the compact case, scattering theory plays a fundamental role in the analysis of $Z(\la)$ and trace formulae when $X$ is not compact but has finite volume, this can be understood by the fact that the Laplacian $\Delta_X$ has continuous spectrum. 
The scattering operator $S(\la)$ is a one-parameter family of linear operators parametrizing the continuous spectrum, this  
is a matrix when ${\rm Vol}(X)<\infty$ and its determinant is a meromorphic function on $\cc$ with a Hadamard 
factorization over its poles and zeros that appears to be of primary importance 
in the study of $Z(\la)$. Indeed, for instance if $\dim X=2$, 
the zeros of $Z(\la)$ in $\Re(\la)<1/2$ are  given by the poles of $\det S(\la)$ with multiplicity
and some poles at $-\nn_0$ of topological order.
Moreover, it can be shown\footnote{See \cite{Mu} for the more general case of finite volume 
surfaces with hyperbolic cusps ends.} that on the critical line $\Re(\la)=1/2$
one can write $\det S(\demi+it)=e^{2\pi i\xi(t)}$ 
for some winding number $\xi(t)\in\rr$ which, added to the counting function of the $L^2$ eigenvalues of $\Delta_X$,
satisfies a Weyl asymptotic formula.\\ 

When $X$ has finite geometry and infinite volume, the continuous spectrum of $\Delta_X$ has infinite multiplicity 
and the scattering operator is not a finite rank operator anymore, making the analysis much more complicated.
Nevertheless, Guillop\'e and Zworski \cite{GZ} did a full and thorough study of that case when $\dim X=2$,
using Birman-Krein theory and the fact that very explicit model operators exist outside a compact set of the hyperbolic surfaces. Similar technics has been used later by Borthwick-Judge-Perry \cite{BJP} to describe the poles and zeros of 
$Z(\la)$ in that case. The higher dimensional case can not be considered directly by this method, for the simple reason that 
there are no good (and natural) model operators near infinity.
However, when $X$ has no cusps ($X$ is called \emph{convex co-compact}), Patterson-Perry \cite{PP} and Bunke-Olbrich \cite{BO} analyzed
the divisors of $Z(\la)$ in $\cc$, while the 
meromorphic extension of $Z(\la)$ follows from dynamical technics of Bowen \cite{B} and Fried \cite{Fr}. 
In that case, $X$ conformally compactifies in a smooth manifold with boundary, the scattering operator 
$S(\la)$ is a pseudo-differential operator on the conformal boundary $\pl\bar{X}$, and
it turns out that the divisors of $Z(\la)$ are given by the poles of $S(\la)$ and the
points in $-\nn_0$ with topological order. In particular, \cite{PP} show
that a zero $\la_0$ of $Z(\la)$ in $\{\Re(\la)<n/2, \la\notin -\nn_0\}$ has order given by
$\tra({\rm Res}_{\la_0}\pl_\la S(\la)S^{-1}(\la))$, where the residue is shown to have finite rank, 
although $\pl_\la S(\la)S^{-1}(\la)$ is definitely not trace class in this case. Considering the result in the finite volume case, one would thus expect $Z(\la)$ to be related to a determinant of $S(\la)$, the sense of which has to be given.\\

One purpose of this work is to carry out such a construction, 
for even-dimensional convex co-compact quotients of $\hh^{n+1}$, 
by developping natural regularization processes that allow to 
understand Selberg trace formula, functional equation and
analysis of the divisors for $Z(\la)$ in a unified way.   
The regularization we need are of two types and are dual in some sense: 
first we regularize integrals of functions in $X$ admitting asymptotic expansions at 
the conformal boundary $\pl\bar{X}$ of $X$, then we regularize trace and determinant of pseudo-differential operators 
on $\pl\bar{X}$ by methods of Kontsevich-Vishik \cite{KV}. 
Interesting facts coming from conformal theory at 
infinity $\pl\bar{X}$ naturally arise in this analysis, and we believe that our approach gives new  
insights about relations between conformal geometry, scattering theory and Selberg zeta function analysis.  
Our method actually allows to deal with much more general geometric settings than convex co-compact 
hyperbolic manifolds, for instance 
we can consider even asymptotically hyperbolic manifolds which are subject 
of recent extensive interest in view of their relations with 
conformal theory and Ads/CFT correspondence.
As a byproduct, we develop a kind of Krein theory through the construction of a 
``counting function of the continuous spectrum'' which we call generalized Krein function.
It is more intrinsic than the usual spectral shift function as we do not need model operators 
to compare. This function can be compared to the winding number $\xi(t)$ described above for finite 
volume manifold, it will be shown to be the phase of the regularized determinant of $S(\la)$ on the continuous spectrum, 
and satisfies Weyl type asymptotics.
Notice that even in the case of Riemann surfaces, our approach is complementary, and probably more intrinsic, 
than that of Guillop\'e-Zworski \cite{GZ} and Borthwick-Judge-Perry \cite{BJP}, at least providing another point of view.\\ 

Let us now give a few definitions before stating the full results. Throughout this paper, $n$ 
will be an odd integer.
An \emph{asymptotically hyperbolic manifold} (AH in short) is defined to be an $(n+1)$-dimensional Riemannian non-compact manifold $(X,g)$
such that $X$ compactifies in a smooth manifold with boundary $\bar{X}=X\cup M$ (here $M=\pl\bar{X}$) and 
there exists a diffeomorphism $\psi:[0,\eps)_x\x M\to \psi([0,\eps)_x\x M)\subset \bar{X}$
such that $\psi|_{x=0}=\textrm{Id}_M$ and
\begin{equation}\label{ahm}
\psi^*g=\frac{dx^2+h(x)}{x^2}
\end{equation} 
for some family of metrics $h(x)$ on $M$ depending smoothly on $x\in[0,\eps)$.
We will say that $g$ is \emph{even modulo} $O(x^{2k+1})$ if the Taylor expansion 
of $h(x)$ at $x=0$ contains only even power of $x$ up to the $x^{2k+1}$ term, this condition
is proved in \cite{G} to be invariant with respect to diffeomorphism $\psi$ and is satisfied 
if the curvature of $g$ is constant outside a compact set of $X$. A choice of $\psi$ 
is actually equivalent to finding a boundary defining function $x$ of $M$ in $\bar{X}$
such that $|dx|_{x^2g}=1$ near $M$ (see \cite{GR,G}), these boundary defining functions will be qualified
of \emph{geodesic}.
A first interesting property of such manifolds is that they induce a conformal class on 
the boundary $M$, the conformal class $[x^2g|_{TM}]$ of $x^2g|_{TM}$ called the \emph{conformal infinity}, 
arising from the non-uniqueness of the boundary defining  
function $x$. Moreover, for any conformal representative $h_0\in[x^2g|_{TM}]$, there exists 
a unique (near $M$) geodesic boundary defining function $x$ of $M$ such that $x^2g|_{TM}=h_0$. 
Finally, a particularly interesting case of AH manifolds is given by those AH metrics 
which satisfies the approximate Einstein equation
$\textrm{Ric}(g)+ng=O(x^{\infty})$, they are nammed \emph{Poincar\'e-Einstein} and
are the main tools of Fefferman-Graham \cite{FGR} theory for studying conformal theory 
of $M$. The main idea is that for a Poincar\'e-Einstein manifold which is even modulo $O(x^{\infty})$, 
all terms $(\pl_x^jh(0))_{j\in\nn}$ are locally determined by $h(0)$ and one can 
extract natural conformal invariants of the conformal infinity of $g$ 
from the studies of Riemannian invariants and natural operators of $g$.\\
  
The Laplacian on an AH manifold has absolutely continuous spectrum 
$\sigma_{\textrm{ac}}(\Delta_g)=[\frac{n^2}{4},\infty)$ and a finite number of 
eigenvalues forming $\sigma_{\textrm{pp}}(\Delta)\subset (0,\frac{n^2}{4})$. 
Mimicking the definition of the counting function of eigenvalues of $\Delta_g$ for compact manifolds, 
one want to define $\xi(t)$ as a trace of the spectral projector
\[\textrm{``}\xi(t)=\tra (\indic_{[\frac{n^2}{4},\frac{n^2}{4}+t^2]}(\Delta_g)) \textrm{''},\]
but this need to be given a precise sense since $\Pi(t^2):=\indic_{[\frac{n^2}{4},\frac{n^2}{4}+t^2]}(\Delta_g)$ is far from being trace class. In Birman-Krein theory \cite{BK,Y}, if one has $2$ positive operators
$P_0, P_1$, the resolvent of which differ by a trace-class operator, 
they induce a distribution $\xi$ defined by 
\[\tra(f(P_1)-f(P_0))=:\int_{0}^\infty \xi(\la)f'(\la)d\la, \quad f\in C_0^\infty(\rr^+)\]
and called \emph{Krein's spectral shift function}, roughly speaking $\xi(t)=\tra(\indic_{[0,t]}(P_1)-\indic_{[0,t]}(P_0))$.
Birman-Krein \cite{BK} show that $-2\pi i\xi$
is the phase of the Fredholm determinant of a 
relative scattering operator of the form $1+K(t)$ with $K(t)$ of trace class on some Hilbert space.
It turns out that Weyl-type asymptotics hold in several standard situations for $\xi$. 
This theory has been shown to be useful
when we want to understand perturbations $P_1$ of a model operator $P_0$, for instance
Guillop\'e and Zworski \cite{GZ} used this theory in an essential way to obtain a 
full scattering analysis of Riemann surfaces of infinite volume, essentially since the geometric decomposition of such manifolds provides model operators near infinity. In our geometric setting, there is definitely a lack of global model near infinity, which is the main difficulty to deduce interesting results from Birman-Krein theory.\\ 

To avoid the problem above, the method we use is to modify the trace functionnal, or more precisely to 
extend it to operators whose integral kernel have certain regularity near infinity. 
To that end, one first needs to extend the notion of integral: if $u$ is a smooth 
function on $X$ which has an asymptotic expansion at the boundary 
\[u(x,y)= \sum_{i=0}^N x^{I+i}u_i(y)+O(x^{\Re(I)+N+1}), \quad u_i\in C^\infty(M), \quad I\in\cc , \quad \Re(I)+N\gg n\]
in a collar neighbourhood $(0,\eps)_x\x M$ of $M$ induced by a geodesic boundary 
defining function $x$ of $M$, then we define its $0$-integral by the following finite part
\[\int^0 u=\textrm{FP}_{s=0}\int x^{s}u\textrm{d}_{g},\]
where the meromorphic extension of $\int x^su\textrm{d}_g$ to $\{\Re(s)>-1\}$ is insured by the expansion 
of $u$ at the boundary. In a fairly natural way, we thus define the $0$-Trace of 
an operators $K$ with continuous Schwartz kernels $\kappa$ by setting 
\[\zerotr(K):=\int^0\kappa|_{\diag}\]
where $\diag$ is the diagonal of $X\x X$, at least when the $0$-integral exists.
Note that $0$-integral and $0$-trace a priori depend on the choice of $x$.\\

Instead of working with the projector $\Pi(t^2)$, we consider its derivative $2t d\Pi(t^2)$ with respect to $t$, which is essentially the spectral measure of $\Delta_g$ on the continuous spectrum, and can be expressed in terms of 
the resolvent by Stone formula.  
The resolvent $R(\la):=(\Delta_g-\la(n-\la))^{-1}$ is bounded on $L^2(X)$ if $\Re(\la)>n/2$, $\la(n-\la)\notin\sigma_{\rm pp}(\Delta_g)$, and Mazzeo-Melrose \cite{MM} analyzed fully 
its integral kernel: in particular they also 
proved\footnote{The condition on evenness has been worked out by \cite{G}}
its meromorphic continuation to $\la\in\cc$ if $g$ is even modulo $O(x^{\infty})$, 
with finite rank polar part at poles, the poles are called \emph{resonances} and their multplicities are defined by
\[m(\la_0)=\left\{\begin{array}{ll}
\rang{\rm Res}_{\la_0} ((2\la-n)R(\la)) & \textrm{ if }\la_0\not=\ndemi\\
\rang{\rm Res}_{\la_0} R(\la) & \textrm{ if }\la_0=\ndemi\end{array}\right.\]
The continuous spectrum is the line $\{\Re(\la)=\ndemi\}$ in this spectral parameter,
and the spectral measure $d\Pi$ is given by Stone formula
\[ d\Pi(t^2)=\frac{i}{2\pi}\Big(R(\ndemi+it)-R(\ndemi-it)\Big), \quad t\in(0,\infty).\]
The precise analysis of $R(\la)$ by \cite{MM} and a refinement to the case of even metrics
allows us prove\footnote{Note that Albin \cite{A1} proved independently that the $0$-Trace of the heat operator on forms renormalizes
with same conditions on the metric. } 
\begin{theo}\label{otrpi}
If $(X,g)$ is asymptotically hyperbolic, even modulo $O(x^{n+1})$, and $n=\dim X-1$ is odd,
then the $0$-Trace of the spectral measure $2td\Pi(t^2)$ is well-defined for $t>0$, extends to $t \in\rr$ 
analytically, to $t\in \cc\setminus (-\nn\cup n+\nn)$ meromorphically 
and is independent of the choice of geodesic boundary defining function $x$.
\end{theo} 

We can thus formally compute $\zerotr(f(\Delta_g))$ for $f\in C_0^\infty(\rr)$ by pairing $f$ with 
the function $t\to \zerotr (2td\Pi(t^2))$, and we naturally call  
\[\xi(t):=\int_{0}^t\zerotr (2ud\Pi(u^2))du \]
the \emph{generalized Krein function} of $\Delta_g$.\\

We want to see its relation with the scattering operator $S(\la)$, in particular we expect,
by comparing with Birman-Krein theory\footnote{Another reason is also given 
in the paper of Carron \cite{CA}, where it is proved that the spectral shift function  
for the comparison of the Laplacian on $X$ with the Laplacian on $X\setminus\{x=\eps\}$ with Dirichlet 
condition is essentially given by the generalized determinant of the Dirichlet-to-Neumann
map on $\{x=\eps\}$.}, that $\xi$ is the phase of a determinant of $S(\la)$ on the continuous spectrum.
In our setting, the scattering operator can be understood as a Dirichlet-to-Neumann
map at infinity, let us recall its definition following \cite{JSB,GRZ}.
For $\Re(\la)=\ndemi$ and $f\in C^\infty(M)$, there exists 
a unique solution of $(\Delta_g-\la(n-\la))u=0$ that satisfies $u=x^\la u_1+
x^{n-\la}u_2$ with $u_i\in C^\infty(\bar{X})$ and $u_2|_{M}=f$. The scattering operator (actually
modified by the Gamma factors) is the map 
\begin{equation}\label{scat}
S(\la): f\in C^\infty(M)\longrightarrow 2^{2\la-n}\frac{\Gamma(\la-\ndemi)}{\Gamma(\ndemi-\la)}u_1|_{M}\in C^\infty(M),
\end{equation}
it is an elliptic pseudo-differential operators (denoted $\psdo$) on $M$ of 
order $2\la-n$, the distributional kernel of which is obtained as weighted 
restriction of the resolvent kernel to $\pl\bar{X}\x\pl\bar{X}$ 
and extends meromorphically in $\la$.\\

Since for $\la>\ndemi$, the  scattering operator $S(\la)$ 
is an elliptic self-adjoint classical $\psdo$ of positive order, with positive 
principal symbol, one can use the method of Kontsevich-Vishik \cite{KV,KV2}, inspired by Ray-Singer's determinant of Laplacian (in \cite{RS}), to define a zeta-regularized determinant of $S(\la)$  by
\[\det S(\la):=e^{-\pl_sZ(\la,s)|_{s=0}},\]
where $Z(\la,s)$ is the meromorphic extension to $\cc$ of the analytic map $s\to \tra(S(\la)^{-s})$ 
defined for $\Re(s)\gg 0$. 

The first main result of this article is   
\begin{theo}\label{det}
If $g$ is even modulo $O(x^{\infty})$ and $n$ odd, the function $\det S(\la)$ extends meromorphically 
to $\cc$ with divisor\footnote{By convention, positive divisors are zeros} at any $\la_0\in\cc$  
\[-m(\la_0)+m(n-\la_0)-\indic_{\ndemi-\nn}(\la_0)\dim \ker S(n-\la_0)+
\indic_{\ndemi+\nn}(\la_0)\dim\ker S(\la_0)\]
Moreover $\det S(\la)$ is a conformal invariant of the conformal infinity $(M,[h_0])$, 
the function $e^{-2i\pi\xi(z)}$ has a meromorphic extension to $\cc$ and
\begin{equation}\label{identdet}
\det S\Big(\ndemi+iz\Big)=(-1)^{m(\ndemi)} e^{-2i\pi\xi(z)}. 
\end{equation}
In particular if $\frac{n^2}{4}-k^2\notin\sigma_{\textrm{pp}}(\Delta_g)$, one has 
\[\det P_k=(-1)^{m(\ndemi)} e^{-2i\pi\xi(-ik)}\]
where $P_k=S(n/2+k)$ is a conformally covariant differential operator on $\pl\bar{X}$, 
which is the $k$-th GJMS conformal Laplacian of \cite{GJMS} when $g$ is Poincar\'e-Einstein. 
\end{theo} 
 
Here the Krein function $\xi$ is expressed as a scattering phase, a phase of $\det S(\la)$,
or a winding number. In the proof, we are actually able to deal with the more general 
case of even metrics modulo $O(x^{n+1})$, the only difference being the domain of meromorphy
given by $\cc\setminus(-\nn\cup n+\nn)$ instead of $\cc$.\\

The main part of the proof is to get a meromorphic identity between the logarithmic derivative of both sides of 
(\ref{identdet}). At least formally, the logarithmic derivative of $\det S(\la)$ with respect to $\la$ 
should be $\pl_\la \log \det S(\la)=\tra(\pl_\la S(\la)S^{-1}(\la))$
but it turns out that $\pl_\la S(\la)S^{-1}(\la)$ is a $\psdo$ of order $\eps,\forall \eps>0$, which is not classical but with
log-polyhomogeneous terms in the local total symbol expansion, and so is far from being trace class. 
To give a sense to the trace of $\pl_\la S(\la)S^{-1}(\la)$, we use the generalized trace functional $\TR$ developped by 
Kontsevich-Vishik \cite{KV,KV2} on classical 
$\psdo$'s on an $n$-th dimensional compact manifold $M$, actually much its extension by Lesch \cite{L} to log-polyhomogeneous
$\psdo$'s, we call $\TR$ the KV-Trace. For a log-polyhomogeneous pseudodifferential operator $A$ on $M$ we set  
$f(A,s):=\tra(AP^{-s})$  where $\Re(s)\gg 0$ and $P$ is any positive self-adjoint elliptic differential operator of order $p\in\nn$, one can then show that $f(A,s)$ has a meromorphic extension in $s\in\cc$ and we set 
\[\TR(A):={\rm FP}_{s=0}f(A,s),\] 
which is a priori dependent of $P$ and is the usual trace if $A$ is trace class. 
We actually prove that $\TR(\pl_\la S(\la)S^{-1}(\la))$ is independent of $P$, meromorphic
in $\la\in \cc$ if $g$ is even modulo $O(x^{\infty})$, 
it is the log-derivative of $\det S(\la)$ when $\la>n/2$, and satisfies the meromorphic identity 
\[-2\pi\pl_z\xi(z)|_{z=i(\ndemi-\la)}=\TR(\pl_\la S(\la)S^{-1}(\la)),\]  
with only simple poles, the residue of which are integers.\\

Applied to even dimensional convex co-compact hyperbolic manifolds, 
this gives another proof of Theorem 1.5 of Patterson-Perry \cite{PP}
relating the divisors of Selberg zeta function to resonances, essentially by showing that $\det S(\la)$ has an explicit 
relation with Selberg zeta function. This relation is through the following functional equation which 
can be compared to the finite volume hyperbolic case \cite[p. 499]{He2}.
\begin{theo}\label{pkhyper} 
Let $X=\Gamma\backslash\hh^{n+1}$ be a convex co-compact quotient of even dimension of $\hh^{n+1}$,
$S_X(\la)$ and $S_{\hh^{n+1}}(\la)$ be the respective scattering operators of $X$ and $\hh^{n+1}$, 
let $P_k$ be the GJMS $k$-th conformal Laplacian of the conformal infinity of $X$, $\chi(\bar{X})$ be 
the Euler characteristic of $\bar{X}$ and $Z(s)$ be the Selberg zeta function of the group $\Gamma$.
Then 
\[\begin{gathered}
\frac{Z(\ndemi-iz)}{Z(\ndemi+iz)}=\frac{\det S_X\Big(\ndemi+iz\Big)}{\Big(\det S_{\hh^{n+1}}\Big(\ndemi+iz\Big)\Big)^{\chi(\bar{X})}},\\
{\rm with }\quad\det S_{\hh^{n+1}}\Big(\ndemi+iz\Big) =\exp\Big(-\frac{2i\pi(-1)^{\frac{n+1}{2}}}{\Gamma (n+1)}\int_{0}^z \frac{\Gamma(\ndemi+it)\Gamma(\ndemi-it)}{\Gamma(it)\Gamma(-it)}dt\Big)
\end{gathered}\]
and if $n^2/4-k^2\notin\sigma_{\textrm{pp}}(\Delta_X)$ with $P_k$ invertible
\[\det P_k=\frac{Z(\ndemi-k)}{Z(\ndemi+k)}\exp\Big(\frac{2\pi(-1)^{\frac{n+3}{2}}}{\Gamma(n+1)}\chi(\bar{X})\int_{0}^k\frac{\Gamma(\ndemi+t)\Gamma(\ndemi-t)}{\Gamma(t)\Gamma(-t)}dt\Big)\]
where the integrals are contour integrals 
avoiding singularities, the final result being independent of the contour. 
\end{theo}

\begin{remark}
It could be noticed that the term $(\det S_{\hh^{n+1}}(\la))^{-\chi(\bar{X})}$ is also the topological 
contribution in the usual functional equation for $Z(\la)$ for compact manifolds (see \cite[3.3.2]{BO1}), 
it did not seem to be remarked in the litterature that this term is a regularized 
determinant of the intertwining Knapp-Stein operator.  
\end{remark}

We conclude with Weyl asymptotic: considering that Weyl type asymptotics are true for spectral shift function, one 
can conjecture  
\[\xi(t)=\frac{(4\pi)^{-\frac{n+1}{2}}}{\Gamma(\frac{n+3}{2})}\zerov(X)t^{n+1}+o(t^{n+1}), \quad t\to +\infty\] 
where $\zerov(X)$ is the renormalized volume defined by the $0$-integral of $1$ (see \cite{GR}).
We actually show in this paper that this holds true (and is even better) in the hyperbolic case, 
\begin{prop}\label{weyla}
For a convex co-compact quotient $X=\Gamma\backslash \hh^{n+1}$ with
dimension of the limit set of $\Gamma$ denoted $\delta$, then we have the Weyl asymptotic
as $t\to \infty$   
\[\xi(t)=
\left\{\begin{array}{ll}
\frac{(4\pi)^{-\frac{n+1}{2}}}{\Gamma(\frac{n+3}{2})}\zerov (X)\Big(t^{n+1} +\sum_{i=1}^{\frac{n-1}{2}}C_it^{2i}\Big)+O(t), & {\rm if }\delta<\ndemi \\
\frac{(4\pi)^{-\frac{n+1}{2}}}{\Gamma(\frac{n+3}{2})}\zerov (X)t^{n+1} +O(t^n), & \rm{ otherwise}
\end{array}\right.\]
where $C_i$ is the $t^{2i}$ coefficient of the polynomial
$\int_{0}^tu\prod_{j=1}^{\frac{n-1}{2}}(\ndemi-j+u^2)du.$ 
\end{prop}
It is worth noticing that a group with a small limit set implies
a much better Weyl asymptotic, we are not aware of such other examples 
in scattering theory.\\

Other natural questions would be 
to deduce an exact trace formula for the wave operator as in \cite{GZ,GN} (for hyperbolic manifolds) 
or \cite{CR} (for asymptotically cylinder manifold), which would require estimates on the counting function for resonances.
Another important step would be to understand the delicate case of $n+1$ odd where things 
do not renormalize correctly.\\  

The paper is organized as follows: we first describe the Kontsevich-Vishik trace
and its extension to odd $\log$-polyhomogeneous $\psdo$'s in odd dimension. Then, 
we recall results of scattering theory, we show the renormalizability of the trace
of the spectral measure and compute this $0$-Trace in function of the 
scattering operator $S(\la)$. We finally define Kontsevich-Vishik determinant 
of $S(\la)$ and its relation with Krein's function.\\

\textbf{Notations}: We use $3$ traces and $2$ determinants in this work, the usual trace (on 
trace class operators) is denoted Tr, the $0$-Trace is written $\zerotr$,
the Kontsevich-Vishik trace is denoted TR, whereas the Fredholm determinant
for operators of the form ``$\textrm{Id}+$ trace class'' is denoted $\det_{\textrm{Fr}}$
and the Kontsevich-Vishik determinant is simply written $\det$.
 
\section{The KV-Trace and odd log-polyhomogeneous pseudo-differential operators}\label{sec2}

\subsection{Log-homogeneous distributions}
A tempered distribution $u$ on $\rr^n$ is said to be log-homogeneous of order $(m,k)\in\cc\x\nn_0$ if for all $t>0$ 
($\mc{S}(\rr^n)$ below is the Schwartz space)
\[\cjg u,\alpha_t^*\psi\cjd=t^{m}\sum_{j=0}^k (\log t)^j\cjg u_j,\psi\cjd, \quad \forall\psi\in\mc{S}(\rr^n)
\]
for some distributions $u_j$ on $\rr^n$ and where $\alpha_t^*\psi(y):=t^{-n}\psi(y/t)$.
Similarly, a smooth function $u$ on $\rr^n\setminus\{0\}$ is said log-homogeneous of order $(m,k)$ if for all $t>0$
we have $u(ty)=t^{m}\sum_{j=0}^k(\log t )^j v_j(y)$
for some smooth $v_j$ on $\rr^n\setminus\{0\}$, this in turn is equivalent to say that $u(y)=\sum_{j=0}^k(\log|y|)^jw_j(y)$ for 
some $w_j(y)$ which are smooth on $\rr^n\setminus\{0\}$ and homogeneous of degree $m$. 
Notice that Fourier transform maps log-homogeneous distributions of order $(m,k)$ to log-homogeneous
distributions of order $(-n-m,k)$.
We say that a log-homogeneous distribution $u$ has negative parity if $\cjg u,\psi(-y)\cjd=-\cjg u,\psi\cjd$ for all 
$\psi\in\mc{S}(\rr^n)$ and positive parity if $\cjg u,\psi(-y)\cjd=\cjg u,\psi\cjd$. 
We use a similar notion for a log-homogeneous function on $\rr^n\setminus\{0\}$: $u$ has negative parity if $u(-y)=-u(y)$, etc...
We first show the following Lemma which extends a classical result in distribution (cf. \cite[Th. 3.2.3]{Ho}).
\begin{lem}\label{distrib}
A log-homogeneous function $u$ on $\rr^n\setminus\{0\}$ of order $(m,k)$ has a unique 
extension as a tempered log-homogeneous distribution of order $(m,k)$ on $\rr^n$ if either $m\not\in -n-\nn_0$, or
if $m\in-n-2\nn_0$ and $u$ has negative parity, or if $m\in-n-1-2\nn_0$ and $u$ has positive parity.  
\end{lem}
\textsl{Proof}: it is similar to Lemma 1.2.5 of \cite{OK}. If $u_z(y):=|y|^zu(y)$ and $\psi\in\mc{S}(\rr^n)$, 
we extend $\cjg u_z(y),\psi\cjd$ from $\Re(z)>-n-m$ to $\Re(z)>-1$ using Taylor expansion of $\psi$ at $0$:  $\psi(y)=\sum_{|\alpha|<N}\pl_y^\alpha \psi(0)y^\alpha/\alpha! +\psi_N(y)$ for $N>-n-m+1$, and computing 
for $\Re(z)\gg 1$:
\[\begin{gathered}
\cjg u_z,\psi\cjd=\int_{|y|>1}|y|^zu(y)\psi(y)dy+\int_{|y|<1}|y|^zu(y)\psi_N(y)dy\\
+\sum_{|\alpha|<N}\sum_{j=0}^k\frac{(-1)^jj!\pl_y^\alpha\psi(0)}{\alpha!(z+m+n+|\alpha|)^{1+j}}\int_{|\theta|=1}
u_j(\theta)\theta^\alpha d\theta
\end{gathered}\]
where $u(y):=|y|^{m}\sum_{j=0}^k(\log|y|)^ju_j(y/|y|)$. The first two terms are holomorphic at $z=0$ and the last term
has a pole at $z=0$ only if $m\in -n-\nn_0$, but then the whole polar part at $z=0$ actually vanishes when $m\in-n-\nn_0$
under our assumptions since, by parity of $u$, $u_j(\theta)\theta^\alpha$ is odd in $\theta$ for $|\alpha|=-m-n$ and for all
$j=0,\dots,k$, so their integral on the sphere vanishes. This proves the holomorphic extension of $u$ since the obtained distribution is log-homogeneous by construction.
As for uniqueness, it suffices to observe that two such distributions would differ from a distribution 
supported at $0$, thus a differential operator applied to the delta function at $0$, which would 
have parity $(-1)^{-m-n}$, thus would be $0$ since $u$ has converse parity by our assumption.
\qed

\subsection{Log-polyhomogeneous pseudo-differential operators}
Let $M$ be a compact manifold of dimension $n$. 
The set $\Psi^m(M)$ ($m\in\rr$) of pseudo-differential operators
on $M$ is defined as follows (see for instance \cite{SH}): 
a continuous operator $A:C^{\infty}(M)\to C^{\infty}(M)$ 
belongs to $\Psi^m(M)$ if its Schwartz kernel $A(y,y')$ is a distribution on $M\x M$, smooth 
outside the diagonal ${\rm diag}_M$
and which can be expressed in coordinates $(y,y')$ near any point $(p,p)\in {\rm diag}_M$ of the diagonal 
under the oscillating integral
\begin{equation}\label{ayy}
A(y,y')=\int e^{i(y-y').\xi}\sigma_A(y,\xi)d\xi
\end{equation}
where $\sigma_A(y,\xi)\in S^m(U)$ is a symbol of order $m$ in a neighbourhood $U$ of $p$, and the 
class $S^l(U)$ consists of the set of $\sigma(y,\xi)\in C^{\infty}(U\x\rr^n)$ satisfying 
\begin{equation}\label{estimeesymbol}
\forall \alpha,\beta \in\nn^n, \exists C_{\alpha,\beta}\geq 0, \quad |\pl_{y}^\alpha\pl_\xi^{\beta}
\sigma(y,\xi)|\leq C_{\alpha,\beta}(1+|\xi|)^{l-|\beta|}.
\end{equation}
An operator $A\in\Psi^l(M)$ is said to be in the class $\Psi^{m,k}(M)$ of log-polyhomogeneous 
$\psdo$'s  of order $(m,k)\in \cc\x\nn_0$ if $\Re(m)<l$ and if its local symbol $\sigma_A(y,\xi)$ 
has an asymptotic expansion (in the usual sense for symbols) of the form 
\begin{equation}\label{totalsym}
\sigma_A(y,\xi)\sim_{|\xi|\to\infty}\sum_{i=0}^{\infty}\sum_{l=0}^ka_{m-i,l}(y,\xi)(\log|\xi|)^l, 
\quad a_{m-i,l}(y,t\xi)=t^{m-i}a_{m-i,l}(y,\xi), \forall t>0.\end{equation}
In \cite{L}, Lesch proves that this condition does not depend on choice of coordinates, 
satisfies the composition law $\Psi^{m,k}(M).\Psi^{m',k'}(M)\subset \Psi^{m+m',k+k'}(M)$,
and also remark that $k=0$ correspond to the classical $\psdo$'s. 
For notational convenience, we will set 
\[a_{m-i}(y,\xi):=\sum_{l=0}^ka_{m-i,l}(y,\xi)(\log|\xi|)^l.\] 

\subsection{Regular parity, odd operators}
We define the class $\Psi_{\reg}^{m,k}(M)$ 
of log-polyhomogeneous $\psdo$'s of order $(m,k)\in\cc\x\nn_0$ 
with regular parity at order $N\in\nn_0$ by the condition that their local total symbol (\ref{totalsym}) satisfies
as functions of $(y,\xi)$
\begin{equation}\label{transmission}
\forall i\leq N,l\leq k,\textrm{ } a_{m-i,l}(y,-\xi)=(-1)^ia_{m-i,l}(y,\xi)\iff 
\forall i\leq N, \textrm{ } a_{m-i}(y,-\xi)=(-1)^ia_{m-i}(y,\xi).\end{equation}
It is straightforward to check that this condition is independent of the choice of coordinates,
indeed the change of symbol under diffeomorphism $\psi$ is given (see Shubin \cite{SH}) by
\[\sigma_{\psi_*A}(\psi(x),\eta)\sim \sum_{\alpha\in\nn_0^n}\frac{1}{\alpha!}(\pl^\alpha_\xi
\sigma_A)(x,\trans D\psi(x)\eta)\Phi_\alpha(x,\eta)\]
where, setting $\psi^{''}_x(z):=\psi(z)-\psi(x)-D\psi(x)(z-x)$,  
$\Phi_\alpha(x,\eta):=\pl_z^\alpha e^{i\psi^{''}_x(z).\eta}|_{z=x}$ which is a polynomial
in $\eta$ of degree less or equal to $|\alpha|/2$. Then writing 
$\sigma_{\psi_*A}\sim\sum_{j\in\nn}b_{m-j}$ with $b_{m-j}$ of order $m-j$ (i.e. including
the log-term) and $\Phi_\alpha(x,\eta)=\sum_{|\beta|\leq |\alpha|/2}c_\beta(x)\eta^\beta$ we get
\[b_{m-j}(\psi(x),\eta)=\sum_{\substack{-|\beta|+i+|\alpha|=j\\
|\beta|\leq |\alpha|/2}}c_\beta(x)\eta^\beta \pl_\xi^\alpha a_{m-i}
(x,\trans D\psi(x).\eta)\]
thus for $j\leq n$, $b_{m-j}(x,-\eta)=(-1)^{j+2|\beta|}b_{m-j}(x,\eta)=(-1)^jb_{m-j}(x,\eta)$ where
we used that $-|\beta|+i+|\alpha|\geq i$. 
The space $\Psi_\reg^{m,k}(M)$ is stable by multiplication by functions in $C^\infty(M)$ and 
composition with operators $\Psi_\reg^{m',k'}(M)$ gives operators in $\Psi_\reg^{m+m',k+k'}(M)$,
indeed it suffices to consider the symbol of the composition of two $\psdo$'s, 
given by (see again \cite{SH}) 
\begin{equation}\label{ab}
\sigma_{AB}(y,\xi)\sim\sum_{j=0}^\infty\sum_{|\alpha|+l+l'=j}\frac{i^{-|\alpha|}}{\alpha!}
\pl^\alpha_\xi a_{m-l}(y,\xi)\pl^\alpha_yb_{m'-l'}(y,\xi)
\end{equation}
and remark that differential operators of even degree have regular parity, thus 
in particular multiplication by smooth functions too.\\

\begin{remark}: For what follows, we shall only need regular parity at order $N=n$, thus from now on we  assume $N=n$.
\end{remark}

A notion of odd classical $\psdo$'s of integer orders ($\geq -n$) has been introduced by Kontsevich-Vishik \cite{KV}
if $M$ has odd dimensional $n$. It corresponds to our definition of regular parity when $k=0$, $N=\infty$ and  $m\in -n+1+2\nn_0$.
Generalizing somehow this notion, we will then say that an operator $A\in\Psi^{m,k}(M)$ with $m\in -n+1+2\nn_0$ 
is odd if it has regular parity at order $n+m$ and $n$ is odd, this defines the odd class $\Psi_\odd^{m,k}(M)$.
Note that $\Psi_\odd^{0,k}(M)=\Psi_\reg^{0,k}(M)$ and that differential operators
of even degree in odd dimension are odd. We will use later the fact that 
if $A\in\Psi_\odd^{m,k}(M)$ and if $B\in\Psi_\odd^{m',0}(M)$ has regular parity at all order,
then $AB\in\Psi_\odd^{m+m',k}(M)$.

Let us now state a useful fact for later 
\begin{lem}\label{power}
If $P\in\Psi_\reg^{d,0}(M)$ with $d>0$, self-adjoint invertible with positive principal symbol, then
$P^{\la}$ and $\log(P)$ can be defined and are respectively in 
$\Psi^{d\la,0}_\reg(M)$ and $\Psi_\odd^{0,1}(M)$. 
\end{lem}
\textsl{Proof}: we proceed exactly as in \cite[Sec. 2]{KV} and \cite[Prop 4.2]{KV}, one can define the power on a cut  
$L_\theta=\{re^{i\theta}, r\in(0,\infty)\}$ for any $\theta\not=0(\pi)$ by 
the method introduced by Seeley (see e.g. Shubin \cite{SH}),
it suffices to write $P^\la$ for $\Re(\la)<0$ as 
\[P^\la=\frac{i}{2\pi}\int_{\Lambda}z^\la (P-z)^{-1}dz\]
where $\Lambda=\{re^{i\theta}, \infty>r\geq\rho\}\cup\{re^{i(\theta-2\pi)},\rho \leq r<\infty)\}
\cup \{\rho e^{i\varphi},\theta>\varphi>\theta-2\pi\}$ for some $\rho>0$ such that the negative eigenvalues
of $P$ are of modulus larger than $\rho$, the power $z^\la$ taken with respect to this cut. 
The polyhomogeneous expansion of $P^\la$ is given (see \cite[Sec. 2]{KV}) by
\[\sigma_{P^\la}\sim\sum_{j=0}^\infty a_{d\la-j}^{(\la)}, \quad 
a^{(\la)}_{d\la-j}(y,\xi)=-(2\pi i)^{-1}\int_\Lambda z^\la b_{-d-j}(y,\xi,z)dz\] 
where $b_{-d-j}(y,\xi,z)$ are homogeneous in $(\xi,|z|^{1/d})$ of degree $-d-j$
and form a complete symbol $\sum_{j}b_{-d-j}$ for $(P-z)^{-1}$.
It is straightforward (for exemple mimicking \cite[Prop. 4.2]{KV}) to see from the construction
of $b_{-d-j}$ that $b_{-d-j}(y,-\xi,z)=(-1)^jb_{-d-j}(y,\xi,z)$ for $j\leq n$ if $P$ has 
regular parity.
Then the homogeneous term $b_{-d-i}$ is transformed into the homogeneous term 
$a_{d\la-j}^{(\la)}(y,\xi)$ with homogeneity $d\la-j$ satisfying the regular parity $a_{d\la-j}^{(\la)}(y,-\xi)=(-1)^ja_{d\la-j}^{(\la)}(y,\xi)$.
The same holds for $\Re(\la)>0$ by mutiplying by some $P^k$ for $k\in\nn$, 
which has regular parity by assumption on $P$. The part with $\log P$
is deduced by differentiating $P^\la$ at $\la=0$, it clearly has regular parity by considering 
equation ($2.11$) of \cite{KV}. 
\qed\\

Remark in this Lemma that we can clearly replace ``regular parity (at order $n$)'' by ``regular parity at order $N$'' for any $N>0$, the proof works as well.

\subsection{Kernels of log-polyhomogeneous operators}
If $A\in \Psi^{m,k}(M)$ with $m\notin -n-\nn_0$ or if $A\in\Psi_{\odd}^{m,k}(M)$ with 
$m\in-n+1-2\nn_0$, the Schwartz kernel\footnote{here we use the notation 
$y:=\pi_L^*y$, $y':=\pi_R^*y$ if $\pi_L,\pi_R:M\x M\to M$ are the left and right projections, 
this notation will be often used along the paper} $A(y,y')$ of $A$  can be decomposed in a neighbourhood $U\x U$ of any point $(y,y)$ of the diagonal into 
\begin{equation}\label{decompnoyau}
A(y,y')=\sum_{i=0}^{N} A_{i}(y,y-y') + A_{N+1}(y,y-y')
\end{equation} 
where $A_{i}(y,u)\in C^{-\infty}(U\x\rr^n)$ are log-homogeneous distributions of order $(-n-m+i,k)$ 
and $A_{N+1}(y,u)\in C^{0}(U\x\rr^n)$. Indeed, 
if $A$ has a symbol $\sigma_A$ with an expansion (as $|\xi|\to \infty$) 
in log-homogeneous functions given by (\ref{totalsym}),  we can write in the distribution sense
$\sigma_A=\sum_{i=0}^Na_{m-i}+(\sigma_A-\sum_{i=0}^Na_{m-i})$ for some large $N>0$ such that 
$(\sigma_A(y,\xi)-\sum_{i=0}^Na_{m-i}(y,\xi))$ is integrable in $\{|\xi|>1\}$, and the $a_{m-i}$ are uniquely determined
log-homogeneous distributions by Lemma \ref{distrib}.
We thus define, using Fourier transform,
\[A_{i}(y,u):=\int e^{iu.\xi}a_{m-i}(y,\xi)d\xi, \quad A_{N+1}(y,u):=\int e^{iu.\xi}\Big(\sigma_A(y,\xi)-\sum_{i=0}^{N}a_{m-i}(y,\xi)\Big)d\xi,\]
the first terms $A_i$ are log-homogeneous distributions in $u$ variable since $a_{m-i}$ are in the $\xi$ variable, 
the $A_{N+1}$ is continuous if $N\geq n+\Re(m)$ since $(1-\chi(\xi))(\sigma_A(y,\xi)-\sum_{i=0}^{N}a_{m-i}(y,\xi))$ is in $L^1(\rr^n,d\xi)$ if $\chi \in C_0^\infty(\rr^n)$ 
equal $1$ near $\xi=0$ and $\chi(\xi)\sum_{i=0}^{N}a_{m-i}(y,\xi)$ has Fourier transform given by convolutions
of log-homogeneous distributions with the Schwartz function $\hat{\chi}$, thus is smooth.
It is also important to notice that, when $A$ has regular parity, $A_i(y,z)$ satisfies the parity rule 
\begin{equation}\label{parityrule}
A_i(y,-z)=(-1)^iA_i(y,z)
\end{equation}
which follows directly from that of the log-homogeneous symbols $a_{m-i}(y,\xi)$.\\

A natural way to consider these kernels is to use polar coordinates $y, r=|y-y'|,\omega=(y-y')/r$ around the diagonal ${\rm diag}_M$ of $M\x M$, this can be formalized globally on $M\x M$ by blowing-up the diagonal. 
Let us recall the blow-up process: if $S$ is a submanifold of a compact manifold $Y$, then consider 
the disjoint union $[Y,S]:=(Y\setminus S)\sqcup SN(S,Y)$ where
$SN(S,Y)\subset TM|_{S}$ is the spherical normal bundle. The blow-down map 
$\beta:[Y,S]\to Y$ is defined to be the identity outside $SN(S,Y)$ and the projection on the
basis on $SN(S,Y)$. The space $[Y,S]$ can be equipped with a structure of smooth manifold with boundary, 
namely the minimal smooth structure for which smooth functions on $Y$ and polar coordinates in $Y$ around $S$
all lift to be smooth \cite[Chap. 5]{ME}. 
For instance, a smooth function on $M\x_0M=[M\x M;\textrm{diag}_M]$ near 
the boundary is the lift by $\beta$ of a 
smooth function on $M\x M\setminus \textrm{diag}_M$ which near $(y,y)\in M\x M$ 
can be written as $f(y,r,w)$ with $f$ smooth and $r=|y-y'|, w=(y-y')/r$; this condition
is independent of the choice of coordinates. Since smooth function $f\in C^\infty(M\x M)$ lifts
under $\beta$ to a smooth function $\beta^*f\in C^\infty(M\x_0 M)$, this 
induces a push-forward $\beta_*:C^{-\infty}(M\x_0M)\to C^{-\infty}(M\x M)$ for distributions defined 
by duality 
$\cjg \beta_*K, f\cjd:=\cjg K,\beta^*f\cjd$
where the pairings are done with respect to a fixed volume density on $M\x M$ and its lift by $\beta$
(or considering in a more invariant way half-densities).  
Then it is clear that an expansion (\ref{decompnoyau}) 
for all $N$ means that the lift $\beta^*(A|_{M\x M\setminus \textrm{diag}_M})$, as a function,
extends to the sum of a function $F\in r^{-m-n}\sum_{j=0}^k(\log r)^jC^{\infty}(M\x_0M)$ and
a function $K\in \beta^*(C^\infty(M\x M))$ where $r$ denotes a global boundary defining function of 
the boundary $SN(\textrm{diag}_M,M\x M)$ in $M\x_0M$. If $\Re(m)>-n$, this gives an $L_{\rm loc}^1$ distribution on 
$M\x_0 M$, the push forward of which is clearly $A$. When $\Re(m)\leq -n$, it is defined as distribution
by holomorphic extension at $z=0$ of the well-defined (in $L_{\rm loc}^1$) distribution $r^z\beta^*A$ 
for $\Re(z)\gg 1$, like in Lemma \ref{distrib}; the extension has no pole at $0$ if $m\notin -n-\nn_0$, 
and it has no pole at $z=0$ either if considered acting on $\beta^*(C^{\infty}(M\x M))$ when $A$ is an odd operator, 
this is easy to check by passing in polar coordinates in the proof of Lemma \ref{distrib}. Moreover 
the push-forward of the obtained distribution under $\beta$ is $A$. This is nothing much more than
reformulating what we said before but set a better ground for our analysis in following chapters.

\subsection{Kontsevich-Vishik trace functional}  
For classical operators $A\in \Psi^{m,0}(M)$ with $m\notin -n-\nn_0$, Kontsevitch and Vishik \cite{KV}
introduced a trace functional, that we denote $\TR$, which extends the usual trace for trace class operators.
They also showed that it keeps a sense when $m\in n-\nn_0$ if $n$ is odd and $A$ is in the odd 
class, as defined above.  
This was generalized later by Lesch \cite{L} for operators in $\Psi^{m,k}(M)$ ($k>0$) when $m\notin -n+\nn_0$.
There are two equivalent ways of defining $\TR(A)$ for $A\in\Psi^{m,k}(M)$ when $m\notin -n-\nn_0$ \cite[Th. 5.6]{L}: 
\begin{itemize} 
\item
The first one is to take any $P\in\Psi^{p,0}(M)$ ($p\in\nn$) 
positive self-adjoint operator with positive principal symbol, 
then extend the function $f(A,s):s\to \tra(AP^{-s})$ meromorphically
from $\Re(s)\gg 1$ to $s\in\cc$, it turns out to be holomorphic at $s=0$ when $m\not\in -n-\nn_0$
and $f(A,0)$ does not depend on $P$, this defines the KV-Trace of $A$ by setting
\begin{equation}\label{defkvtrace}
\TR(A):=f(A,0).
\end{equation}
\item The second way is to define the ``density'' (using notations (\ref{totalsym}) for $\sigma_A$)  
\begin{equation}\label{omegakv}
\omega_{\textrm{KV}}(A)(y):=(2\pi)^{-n}\Big(\textrm{FP}_{\eps\to 0} \int_{|\xi|<\eps^{-1}} \sigma_A(y,\xi)d\xi\Big) |dy|.
\end{equation}
which turns out to be a true density, i.e. independent on choice of coordinates, and to set 
\[\TR(A):=\int_M\omega_{\rm KV}(A)\]
\end{itemize}
Note that in both case $\TR$ is the usual trace if $\Re(m)<-n$.
When $m\in-n+\nn_0$, Lesch shows, extending \cite{KV} when $k>0$, that 
$f(A,s)$ extends meromorphically with a possible pole of order $\leq k+1$ 
whose $(k+1)$-th coefficient in the polar part of Laurent expansion at $0$ is
the $k$-th Wodzicki residue \cite[Cor. 4.8, Th. 5.6]{L}
\begin{equation}\label{wod}
\textrm{WRes}_k(A)=\frac{(k+1)!}{(2\pi)^n}\int_{S^*M}a_{-n,k}(y,\xi)|dyd\xi|\end{equation}
with notations of (\ref{totalsym}) where $S^*M$ is the unit bundle of $T^*M$, 
this number is proved to be globally well-defined. 
\begin{lem}\label{indept}
Suppose that $A\in\Psi_{\odd}^{m,k}(M)$ with $m\in-n+1+2\nn_0$ or $A\in\Psi^{m,k}(M)$ with $m\notin -n+\nn_0$.
Then $\omega_{\rm KV}(A)$ defined by \eqref{omegakv} is a density and the function $s\to \TR(AP^{-s})$
is holomorphic at $0$ if $P\in\Psi^{p,0}(M)$ has 
regular parity at order $n+m$. The KV-Trace of $A$ can then be defined equivalently 
by $\TR(AP^{-s})|_{s=0}$ or by $\TR(A)=\int_M\omega_{\rm KV}(A)$, 
the result is independent of $P$, linear in $A$. 
It is a trace in the sense that 
\begin{equation}\label{traceab}
\TR(AB)=\TR(BA), \quad AB, BA\in(\cup_{m\in \cc\setminus -n+\nn_0}\Psi^{m,k}(M))\cup (\cup_{m\in-n+1+2\nn_0}\Psi_\odd^{m,k}(M)).\end{equation}
and is the usual trace on trace class operators, i.e on $\Psi^{m,k}(M)$ if $m<-n$. 
Finally, if $A=A(\la)\in\Psi_\odd^{m,k}(M)$ depends analytically on a parameter $\la\in\cc$ in the 
sense of \cite[Def. 1.9]{PS}, with order $(m,k)$ constant in $\la$, then $\TR(A(\la))$ is analytic in $\la$.
\end{lem}
\textsl{Proof}: we use notation (\ref{totalsym}) and set $m\in-n+1+2\nn_0$ since the other cases are proved in 
\cite[Th. 5.6]{L}. Actually the proof is also along the lines of the paper of Lesch \cite{L}. 
The term $\omega_{KV}(A)$ is independent of coordinates, this is obtained exactly like 4) of Lemma 5.3 in \cite{L}
using the linear change of coordinates computed in Proposition 5.2 of \cite{L}: indeed the defect (to be a density) in 
the change of coordinates vanishes when $A$ has regular parity to the right order since it involves only the integrals 
of $a_{-n,l}(x,\xi)$ times even functions of $\xi$ on $\{\xi=1\}$, which are thus odd functions for all $l\leq k$.
We have for $\Re(s)\gg 0$
\begin{equation}\label{traps}
\TR(AP^{-s})=\tra(AP^{-s})=\int_M\omega_{\textrm{KV}}(AP^{-s}),\end{equation}
the last identity holds obviously (see \cite[Eq. 5.19]{L}) since the operator 
is trace class. We have to prove that $\omega_{\textrm{KV}}(AP^{-s})$ extends holomorphically 
to $s=0$, this is actually quite straightforward using the proof of \cite[Lem. 5.4]{L} and 
the vanishing of the integral $a_{-n,l}(x,\xi)$ on
$\{|\xi|=1\}$, and $\TR(A)=\int_M\lim_{s=0}\omega_{\textrm{KV}}(AP^{-s})$.
Indeed, \cite[Lem. 5.4]{L} shows that the density $\omega_{\textrm{KV}}(AP^{-s})$ is equal to
\begin{equation}\label{wkvaps}
\int_{\rr^n}a^{(s)}_N(y,\xi)d\xi+\sum_{j=0}^N\sum_{l=0}^k\int_{|\xi|<1}\psi(|\xi|)a^{(s)}_{z(s)-j,l}(y,\xi)d\xi+
\frac{(-1)^ll!}{(z(s)+n-j)^{l+1}}\int_{|\xi|=1}a^{(s)}_{z(s)-j,l}(y,\xi)d\xi
\end{equation}
times $(2\pi)^{-n}|dy|$ if $AP^{-s}$ has the symbol expansion $\sum_{j=0}^{N}\sum_{l=0}^ka^{(s)}_{z(s)-j,l}(y,\xi)(\log|\xi|)^l+a^{(s)}_{N}(y,\xi)$ with
$a^{(s)}_N=O(|\xi|^{-n-1})$ as $|\xi|\to\infty$, $z(s):=-ps+m$ with all terms holomophic 
in $s$ near $0$ 
and $\psi$ a cutoff with support in $[1/2,\infty]$ which equals $1$ in $[1,\infty)$. 
In (\ref{wkvaps}), the terms with $j=n+m$ in the integral on $\{|\xi|=1\}$ 
turn out to be $0$, this is a consequence 
of the fact that $AP^{-s}\in\Psi^{z(s),0}(M)$ has regular parity at order $n+m$ near $s=0$ 
(thus $a^{(s)}_{z(s)-n-m,l}$ is odd in $\xi$, $n+m$ being odd) by assumption on $P$, the remark 
following Lemma \ref{power} and the multiplicative property of regular parity at order $n+m$. 
This proves the holomorphic extension at $s=0$ of $\omega_{\textrm{KV}}(AP^{-s})$ and
the independence with respect to $P$ since the value of (\ref{wkvaps}) at $s=0$ 
depends only on $a^{(0)}_N$ and $a^{(0)}_{m-j,l}$. This also gives 
$\lim_{s=0}\omega_{\textrm{KV}}(AP^{-s})=\omega_{\textrm{KV}}(A)$.
The fact about the cyclicity of $\TR$ 
is proved in \cite[Prop. 3.2]{KV} for classical operators such that the order of $AB$
is not integer, their proof applies word by word in our case too, using
the expression of the KV-density given in the next Lemma \ref{kvkernel}. 
As for the last statement about analyticity, we can introduce an analytic dependence 
of $A$ (in the $\psdo$ topology, see e.g. \cite[Def. 1.9]{PS} for definition) with respect to some paramater $\la\in\cc$ so that the order of $A$ remains constant
with respect to $\la$ and that $A$ is an operator for all $\la$, then the expression (\ref{wkvaps}) taken at $s=0$
is clearly analytic in $\la$. 
\qed\\

To compute the KV-Trace of $A\in\Psi_\odd^{-n+m,1}(M)$ with $m\in 1+2\nn_0$, it suffices by linearity to compute 
the KV-Trace of operators supported in charts of the manifold, or in other words
to consider operators $A\in\Psi_\odd^{-n+m,1}(U)$ in a bounded open set 
$U\subset \rr^n$, with compactly supported kernel in $U\x U$. We search to express the 
KV-Trace in term of the Schwartz kernel, in a way similar to \cite[Lem. 2.2.1]{OK}.
\begin{lem}\label{kvkernel}
Let $U\subset \rr^n$ be an open set, let $A\in\Psi_\odd^{m,k}(U)$ 
with Schwartz kernel $A(y,y')$ compactly supported in $U\x U$. Set    
$A_i \in C^{-\infty}(U\x\rr^n)$ for $i=0,1,\dots,n+m$  and
$A_{n+m+1}\in C^{0}(U\x\rr^n)$ defined by \eqref{decompnoyau},
then the KV density of $A$ is given by
\begin{equation}\label{traomega}
\omega_{\rm{KV}}(A)(y)=A_{n+m+1}(y,0)|dy|
\end{equation}
\end{lem}
\textsl{Proof}: let us consider $m=0$ for simplicity and since it will 
be the case of interest later, the other cases are obviously similar. 
Let $\psi\in C_0^\infty([0,\eps))$ which equals $1$ near $0$ and so that $\psi(|y-y'|)A(y,y')=A(y,y')$, 
then $(y.y')\to \psi(|y-y'|)A_{n+m+1}(y,y-y')$ is the kernel of a trace class
$\psdo$, thus its KV density is the restriction on the diagonal $A_{n+m+1}(y,0)$.  
Now it suffices then to prove that the distribution 
\[ (y,y')\to \psi(|y-y'|)\sum_{i=0}^{n}A_i(y,y-y')\]
has vanishing KV density. The symbol of this operator is, up to constant, given by the Fourier transform
$\sigma(y,\xi):=\mc{F}_{z\to-\xi}(\psi(|z|)\sum_{i=0}^nA_i(y,z))$. We decompose this distribution into
\[\sigma(y,\xi)=\sum_{i=0}^n\mc{F}_{z\to -\xi}A_{i}(y,z) + \mc{F}_{z\to -\xi}((\psi(|z|)-1)A_i(y,z)).\]
Let us denote by $a_{-i}(y,\xi)$ the first term in the sum and by $\beta_{-i}(y,\xi)$ the second.
Since $\sigma(y,\xi)$ is smooth, the density $\omega_{\rm KV}(A)(y)={\rm FP}_{R\to\infty}\int_{|\xi|<R}\sigma(y,\xi)d\xi$ 
can also be defined by 
\[\omega_{\rm KV}(A)(y)={\rm FP}_{R\to\infty}\int_{1/R<|\xi|<R}\sigma(y,\xi)d\xi.\] 
The $a_{-i}(y,\xi)$ is log-homogeneous in $\xi$ of degree $(-i,k)$, thus $L^1_{\rm loc}$ if $i<n$, and can be written
under the form $a_{-i}(y,\xi)=\sum_{l=0}^ka_{-i,l}(y,\xi)(\log |\xi|)^l$ for some $a_{-i,l}(y,\xi)$ homogeneous
of degree $-i$ in $\xi$. A straightforward computation yields for $i<n$
\[\int_{1/R<|\xi|<R}a_{-i}(y,\xi)d\xi=\sum_{l=0}^k\pl^l_s
\left(\frac{R^{s-i+n}-R^{-s+i-n}}{s-i+n}\right)|_{s=0}\int_{S^n}a_{-i,l}(y,\theta)d\theta ,\]
which clearly has vanishing finite part as $R\to 0$. As for $b_{-i}(y,\xi)$ when $i<n$, it is the Fourier transform of
a log-polyhomogeneous symbol of order $(-n+i,k)$, so it is in $L^1$ in $\xi$ and $\int_{1/R<|\xi|<R}b_{-i}(y,\xi)d\xi$
has a limit (its integral on $\rr^n$) as $R\to \infty$, which in turn is the value of $(\psi(|z|)-1)A_i(y,z)$ at $z=0$, that is $0$.

Let us finally consider the case $i=n$.  We have $A_n(y,-z)=(-1)^nA_n(y,z)$ by (\ref{parityrule}), so $\psi(|z|)A_n(y,-z)=
(-1)^n\psi(|z|)A_n(y,z)$ and its Fourier transform $z\to \xi$ has the same parity law in $\xi$. Then
\[\int_{1/R<|\xi|<R}(a_{-n}+b_{-n})(y,\xi))d\xi=\int_{1/R}^R\int_{S^n}(a_{-n}+b_{-n})(y,r\theta)d\theta r^{n-1}dr\]
and a change of variable $\theta\to -\theta$ in the integral together with the parity
law of of $a_{-n}+b_{-n}$ shows that this integral vanishes for all $R>1$. This ends the proof.
\qed\\
 
We apply this result to the following
\begin{cor}\label{corol}
Let $U\subset \rr^n$ an open set, let $A\in\Psi_\odd^{0,k}(U)$ with compactly supported kernel $A(y,y')$ in $U\x U$. Assume
that there exists some functions $W\in C^{\infty}(U)$ and some log-homogeneous $B_i \in C^{\infty}(U,\rr^n\setminus\{0\})$,
$i=0,\dots,n$, of order $(-n+i,k)$ satisfying $B_{i}(y,-z)=(-1)^iB_{i}(y,z)$
and such that
\[\Big|A(y,y+z)-W(y)-\sum_{i=0}^{n}r^{-n+i}B_i(y,z)\Big|=O(|z|^{\demi}) \quad |z|\to 0.\] 
Then $\omega_{\textrm{KV}}(A)= W|dy|$.
\end{cor}
\textsl{Proof}: If $A_i$ are defined by \eqref{decompnoyau}, then asymptotic considerations near $z=0$ 
show that, as functions on 
$U\x (\rr^{n}\setminus 0)$, $A_{i}=B_i$ for $i<n$, and finally 
\[(A_{n,0}(y,z)-B_{n,0}(y,z))|dy|=W(y)|dy|-\omega_{\textrm{KV}}(A)(y), \quad \forall z \textrm{ near }0,\]
where $A_{n}(y,z)=\sum_{l=0}^kA_{n,l}(y,z)(\log|z|)^l$ for some $A_{n,l}(y,z)$ homogeneous in $z$ of 
degree $0$, and the obvious similar decomposition for $B_n$.
In particular, since $A_n(y,z),B_n(y,z)$ are odd in $z$, 
it is immediate to see that $W|dy|=\omega_{\textrm{KV}}(A)$.
\qed

\subsection{KV-Trace of a product}
In this subsection, we express the KV trace of a product of $2$ pseudodifferential operators.
We state the result and postpone the proof, which is a bit technical, to the Appendix. 

We consider a case which will be of special interest later, that is the composition of $2$
polyhomogeneous pseudo-differential operators $F\in\Psi_{\reg}^{2\la-n,1}(M)$ and $L\in\Psi_{\reg}^{n-2\la,0}(M)$
for $\la\not\in\demi\zz$, then the composition $FL\in\Psi_{\odd}^{0,1}(M)$ has a well-defined KV-Trace. In our application,
$F$ will be the derivative $\pl_\la S(\la)$ of the scattering operator while $L$ will be its inverse 
$S^{-1}(\la)$, keeping in mind that we wish to study $\TR(\pl_\la S(\la)S^{-1}(\la))$, see next section.\\

Let us first introduce a notation to simplify statements.
\begin{defi}\label{regusingu}
If $U\subset \rr^n$ is a bounded open set and $u\in C^{\infty}(U\x(0,1)\x S^{n-1})$
has an asymptotic expansion as $r\to 0$
\[u(y',r,w)\sim \sum_{i=0}^{\infty}\sum_{l=0}^k r^{-n-1+i}(\log r)^{l}u_{i,l}(y',w), \quad 
u_{i,l}\in C^\infty(U\x S^{n-1}),\]
we define 
\begin{equation}\label{usingreg}
[u]_{\sing}:=\sum_{i=0}^{n-1}\sum_{l=0}^k r^{-n-1+i}(\log r)^{l}u_{i,l},\quad 
[u]_{\reg}:=u-\sum_{i=0}^{n}\sum_{l=0}^k r^{-n-1+i}(\log r)^{l}u_{i,l}.\end{equation}
\end{defi}

For smoothing pseudodifferential operators $F,L \in \Psi^{-\infty}(M)$, it is well known 
that the usual trace is the integral of the kernel of $FL$
on the diagonal
\[\tra(FL)=\int_{M\x M} F(y,y')L(y',y)\textrm{d}_{h_0}(y')\textrm{d}_{h_0}(y).\]
For our case, this is singularly different and the whole problem comes from the diagonal singularities
of the kernels. One way to approach it is too work in local coordinates near the diagonal and to blow
it up, i.e. to use polar coordinate around the diagonal.  

\begin{prop}\label{kvtrace}
Let $(M,h_0)$ be a compact Riemannian manifold , $F\in\Psi_{\reg}^{2\la-n,1}(M)$ and $L\in \Psi_{\reg}^{n-2\la,0}(M)$ with symmetric Schwartz
kernel $F(y,y')=F(y',y)$ and $L(y,y')=L(y',y)$ outside the diagonal $\{y=y'\}$.
Let $\mc{U}$ be an atlas and $(\chi_j)_{j\in J}$ an associated partition of unity of $M$ such that  
if $i,j$ satisfy $\supp\chi_i\cap \supp\chi_j\not=\emptyset$, there exists a chart $U_{ij}$
that contains both supports.
Then using Definition \ref{usingreg}, there exists $A>0$ 
such that the KV-Trace of $FL$ is, for any $A'>A$,
\begin{equation*}
\begin{gathered}
\TR(FL)=\sum_{\substack{i,j \in J\\
\supp\chi_i\cap\supp\chi_j\not=\emptyset}}\Big(\int_{U_{ij}}\int_{0}^{A'}\int_{S^{n-1}}\Big[\beta^*(\chi_i|\det h_0|^{\demi}\otimes\chi_j)\beta^*F
\beta^*L\Big]_{\reg}(y,r,w)dwdr\textrm{d}_{h_0}(y)\\
\quad\quad\quad\quad\quad\quad-\int_{U_{ij}}\int_{A'}^\infty\int_{S^{n-1}}\Big[\beta^*(\chi_i|\det h_0|^{\demi}\otimes\chi_j)\beta^*F
\beta^*L\Big]_{\sing}(y,r,w)dwdr\textrm{d}_{h_0}(y)\Big)\\
+ \sum_{\substack{i,j \in J\\
\supp\chi_i\cap\supp\chi_j=\emptyset}}\int_{U_{ij}} \chi_i(y)\chi_j(y')F(y,y')L(y',y) \textrm{d}_{h_0}(y)\textrm{d}_{h_0}(y').
\end{gathered}
\end{equation*}
where $\beta$ the blow-down map $\beta: (y,r,w)\to (y,y')=(y,y+r\omega)$ mapping 
$U_{ij}\x [0,\infty)\x S^{n-1}$ into $U_{ij}\x \rr^n$ .
\end{prop}

\section{Krein's formula}

\subsection{Renormalized integral}\label{renormalized}
Let us recall a result of Graham-Lee \cite{GRL} which says
that for any $h_0\in [x^2g|_{TM}]$, there exists a unique (in a neighbourhood of $M$) 
smooth boundary defining function $x$ of $M=\pl\bar{X}$ in $\bar{X}$ such that the metric 
has the form, in a collar neighbourhood $(0,\eps)_x\x M$ of $M$,
\begin{equation}\label{modelform}
g=\frac{dx^2+h(x)}{x^2}, \quad h(0)=h_0\end{equation}
for $h(x)$ a smooth $1$-parameter family of metrics on $M$. Then $x$ will be called the 
\emph{geodesic boundary defining function} correponding to $h_0$.\\

We now recall a couple of things on renormalized integrals. Let $x$ be a geodesic boundary defining function
of $\bar{X}$. As written in the introduction, if $u\in C^{\infty}(X)$ has an expansion in a collar neighbourhood 
$(0,\eps)_x\x M_y$ of $M$ 
\[u(x,y)=\sum_{i=0}^Nx^{I+i}u_i(y)+O(x^{\Re(I)+N+1}), \quad I\in\cc\]
for some $N>-\Re(I)$, we define the $0$-integral of $u$ by
\[\int^0u:=\textrm{FP}_{t\to 0}\int_{x>t}u(x,y)\textrm{d}_{g}\]
where $\textrm{FP}_{t\to 0}$ means ``finite part at $t=0$'', 
this is the constant term in the expansion, after noticing that
there is always an expansion in powers of $t$ and $\log t$ 
with a remaining term being $O(t)$ by assumption on $u$. The $0$-integral depends
a priori on the function $x$.
Albin \cite[Sec. 2.2]{A} proved that if $I=0$ then
\[\int^0u=\textrm{FP}_{z=0}\int_Xx^{z}u\textrm{d}_g, \quad z\in \cc,\]
here the finite part is the regular part at $z=0$ in the Laurent expansion at $z=0$,  
a straightforward computation shows that it also holds if $I\in\cc$.
Actually, if $I\notin n-\nn_0$, then the $0$-integral of $u$ is independent of the choice 
of boundary defining function $x$. Indeed, following \cite[Prop. 2.2.]{A}, if $\hat{x}=e^\omega x$
is another geodesic boundary defining function in $\bar{X}$ then
\[\int_X (x^z-\hat{x}^z)u \textrm{d}_g=z\int_Xx^z\frac{1-e^{\omega z}}{z}u\textrm{d}_g\]
and the integral on the right is holomorphic at $z=0$, this is easily checked
by using an expansion of the integrand in (non-integer) powers of $x$, thus the right hand side
vanishes at $z=0$. If $I\in n-\nn_0$ the $0$-integral depends on $x$ but it is proved in
\cite[Th. 2.5]{A} that for $n$ odd and $I=0$, if $h(x)$ in (\ref{modelform}) has an even 
expansion in $x$ modulo $O(x^{n})$ with $\tra_{h_0}(\pl^n_xh(0))=0$ and if $u$ has
an even expansion in $x$ modulo $O(x^{n+1})$, then the $0$-integral of $u$ is independent of the choice
of $x$.
   
\subsection{Resolvent} \label{resolvent}

We first recall definitions of blow-ups $\bar{X}\x_0\bar{X}$, $\bar{X}\x_0 M$, $M\x_0M$,
as defined for instance in \cite{MM,JSB} and the class of pseudo-differential
operators on $X$ which contains the resolvent of the Laplacian. 
We will generally denote by $[Y;S]$ the (normal) blow-up of $Y$ around a submanifold $S$, 
the spherical normal interior pointing bundle of $S$ in $Y$ is called front face of the blow-up
and the canonical map $\beta:[Y;S]\to Y$ which is the identity outside the front face and the projection
on the base on the front face is called blow-down map. 
In first section we introduced briefly this concept for a submanifold $S$ embedded 
in a compact manifold $M$, it can be actually generalized when $S$ is a submanifold
with corners of a manifold with corners.
For our case, the blow-up $\bar{X}\x_0\bar{X}$ is defined as a set
\[\bar{X}\x_0\bar{X}:=[\bar{X}\x\bar{X};\textrm{diag}_M]=(\bar{X}\x\bar{X}\setminus \textrm{diag}_M)\sqcup
SN_+(\textrm{diag}_M,\bar{X}\x\bar{X})\]
where $\textrm{diag}_M$ is the diagonal in $M\x M\subset \bar{X}\x\bar{X}$, 
$SN_+(\textrm{diag}_M,\bar{X}\x\bar{X})$ is the spherical normal interior pointing
bundle of $\textrm{diag}_M$ in $\bar{X}\x \bar{X}$, the front face of $\bar{X}\x_0\bar{X}$ is denoted $\mc{F}$,
the blow-down map is denoted by 
\begin{equation}\label{beta}
\beta:\bar{X}\x_0\bar{X}\to \bar{X}\x\bar{X}.\end{equation}
A topological and smooth structure of manifold with corners can be given on $\bar{X}\x_0\bar{X}$
through normal fibrations of $\textrm{diag}_M$ in $\bar{X}\x\bar{X}$
and polar coordinates: for instance if $(x,y)$ (here $y=(y_1,\dots,y_n)$) are coordinates on a 
neighbourhood of $y_0\in M$, then defining $(x,y):=\pi^*_L(x,y)$ and $(x',y'):=\pi_R^*(x,y)$ with 
\begin{equation}\label{pirpil}
\pi_L,\pi_R: \bar{X}\x\bar{X}\to \bar{X}, \quad \pi_R(m,m')=m', \quad \pi_L(m,m')=m\end{equation}
gives coordinates $(x,y,x',y')$ near $(y_0,y_0)\in \textrm{diag}_M$ and a function $f$  
on $\bar{X}\x_0\bar{X}$ supported near the fibre $\mc{F}_p$ of $\mc{F}$ where $p=(y_0,y_0)\in\textrm{diag}_M$
is said smooth if $\beta_*f$, defined outside $\textrm{diag}_M$, can be expressed as a smooth function of the 
polar variables
\begin{equation}\label{coordeclate}
R=(x^2+{x'}^2+|y-y'|^2)^{\demi}, \quad\rho:=\frac{x}{R}, \quad \rho'=\frac{x'}{R}, \quad 
\omega:=\frac{y-y'}{R}, \quad y'.
\end{equation}
This manifold with corners has three kind of boundary hypersurfaces, the front face $\mc{F}$
defined in these coordinates by $R=0$, it is a bundle in quarter of sphere, the two other 
boundary faces are called top and bottom faces, in these coordinates $\mc{T}:=\{\rho=0\}, \mc{B}:=\{\rho'=0\}$.
The diagonal of $X\x X$ lifts under $\beta$ to a submanifold, denoted $\textrm{diag}$, whose closure 
is denoted $\pl\textrm{diag}$ and meets only the topological boundary of $\bar{X}\x_0\bar{X}$ at $\mc{F}$. Coordinates
(\ref{coordeclate}) are actually not really coordinates in the usual sense, one actually has to
consider three systems of coordinates induced by $(x,y,x',y')$ that cover $\mc{F}$ near 
$\mc{F}_p$. The first two ones are projective coordinates
\[s:=\beta^*\Big(\frac{x}{x'}\Big), \quad z:=\beta^*\Big(\frac{y-y'}{x'}\Big), \quad x':=\beta^*(x'), \quad 
y':=\beta^*(y')\]
\[t:=\beta^*\Big(\frac{x'}{x}\Big), \quad z':=\beta^*\Big(\frac{y-y'}{x}\Big), \quad x:=\beta^*x, \quad 
y:=\beta^*y\]
covering (near $\mc{F}_p$) respectively a neighbourhood of $\mc{F}\cap\mc{T}, \mc{F}\cap\mc{B}$ and 
valid on the whole interior of the front face $\mc{F}$ (near $\mc{F}_p$). Note that $\mc{F}$ is repectively
$\{x'=0\},\{x=0\}$ in these coordinates and $\textrm{diag}=\{s=1,z=0\},\{t=1,z'=0\}$. The last system is 
\[w:=\beta^*\Big(\frac{y-y'}{r}\Big), \quad \varrho:=\beta^*\Big(\frac{x}{r}\Big),\quad
\varrho':=\beta^*\Big(\frac{x'}{r}\Big), \quad r:=\beta^*(|y-y'|), 
\quad y':=\beta^*(y'),\]
covers a neighbourhood of $\mc{T}\cap\mc{B}\cap\mc{F}$ but is not defined at $\pl\textrm{diag}$.
If $R$ is a global defining function of $\mc{F}$ in $\bar{X}\x_0\bar{X}$, we have canonically associated
global boundary defining functions $\rho:=x/R,\rho'=x'/R$ of respectively $\mc{T},\mc{B}$.
The blow-up $\bar{X}\x_0M=[\bar{X}\x M;\textrm{diag}_M]$ is defined similarly 
and it is direct to see that it is canonically diffeomorphic to the bottom face $\mc{B}$ 
of $\bar{X}\x_0\bar{X}$. Its front face is denoted $\mc{F}'\simeq \mc{F}\cap\mc{B}$ and the blow-down map
$\beta'$. 
Then the final blow-up $M\x_0M:=[M\x M;\textrm{diag}_M]$  
is canonically diffeomorphic to $\mc{T}\cap\mc{B}$, its front face is denoted $\mc{F}_\pl\simeq 
\mc{F}\cap\mc{T}\cap\mc{B}$ and the blow-down map $\beta_\pl$.\\

The space of smooth functions vanishing at all order on a manifold with corners $Y$ is 
denoted by $\dot{C}^{\infty}(Y)$, its topological dual with respect to a pairing
induced by a given volume density is the space of extendible distributions $C^{-\infty}(Y)$.
Since $\beta^*:\dot{C}^\infty(\bar{X}\x\bar{X})\to \dot{C}^{\infty}(\bar{X}\x_0\bar{X})$
is an isomorphism, there is an induced isomorphism $\beta^*$ between their duals, the push-forward 
$\beta_*$ is its inverse isomorphism between $C^{-\infty}(\bar{X}\x_0\bar{X})$ and
$C^{-\infty}(\bar{X}\x\bar{X})$. Schwartz Theorem \cite{ME1} asserts
in this setting that the set of continuous operators from $\dot{C}^{\infty}(Y)$
to $C^{-\infty}(Y')$ is in one to one correspondence with $C^{-\infty}(Y\x Y')$ if
$Y,Y'$ are smooth manifolds with corners and assuming volume densities are given 
(otherwise introduce half-densities).\\
 
We define, following \cite{MM}, a natural class of differential operators
on $\bar{X}$ degenerating uniformly at $M$: the set 
$\textrm{Diff}_0^m(\bar{X})$ for $m\in\nn_0$ is the set of smooth differential operators
on $\bar{X}$ of order $m$ such that $P\in\textrm{Diff}^m_0(\bar{X})$ if it can be written as
\[P=\sum_{j+|\beta|\leq m}P_{j,\alpha}(x,y)(x\pl_x)^j(x\pl_y)^\beta, \quad P_{j,\alpha}\in C^\infty(\bar{X})\] 
locally near the boundary $M=\{x=0\}$.
The Laplacian on an asymptoticaly hyperbolic manifold $(X,g)$ is an operator in $\textrm{Diff}_0^2(\bar{X})$,
it can be expressed in a collar neighbourhood arising from a geodesic boundary defining function $x$ by
\begin{equation}\label{laplacian}
\Delta_g=-(x\pl_x)^2+nx\pl_x+x^2\Delta_{h(x)}-\demi x\tra_{h(x)}(\pl_xh(x))x\pl_x.
\end{equation}

Let $R,\rho,\rho'$ be boundary defining functions of the three boundary hypersurfaces $\mc{F},\mc{T},\mc{B}$.
We define $R^j\Psi_0^{m,k,l}(\bar{X})$ for $m,k,l\in\cc, j\in\nn_0$ to be the set of continuous linear operators $A$ 
from $\dot{C}^{\infty}(\bar{X})$ to its dual $C^{-\infty}(\bar{X})$ (pairing through
the volume density of $g$) such that the lift of the Schwartz kernel $\kappa_A$ of $A$ 
by $\beta$ satisfies 
\[\beta^*\kappa_A\in \rho^{k}{\rho'}^lR^{j}C^{\infty}(\bar{X}\x_0\bar{X})+ R^jI^{m}(\bar{X}\x_0\bar{X},\textrm{diag})\] 
where $I^m(\bar{X}\x_0\bar{X},\textrm{diag})$ stands for the set of distributions classically conormal of order $m$ 
to the interior diagonal, i.e. those which can be written locally as Fourier transform of symbols of order $m$
in the fibres of the normal bundle $N(\textrm{diag},\bar{X}\x_0\bar{X})$ of $\textrm{diag}$ through a normal fibration.
We will also denote by $\Psi^{k,l}(\bar{X})$ the operators having a Schwartz kernel in $x^k{x'}^lC^{\infty}
(\bar{X}\x\bar{X})$.\\   

We know from \cite{MM,G} that the modified resolvent $R(\la)=(\Delta_g-\la(n-\la))^{-1}$ is a bounded operator on $L^2(X)$
for $\Re(\la)>\ndemi$ and $\la(n-\la)\notin \sigma_{pp}(\Delta_g)$ which 
extends to $\la\in\cc\setminus \{(n-1)/2-k-\nn\}$ meromorphically with poles of finite multiplicity (the rank of
the polar part in Laurent expansion at a pole is finite) if the metric is even modulo $O(x^{2k+1})$, this results holds for any $k\in\nn$. The resonances are the poles of $R(\la)$ with finite multiplicity, the multiplicity is defined by
\[m(\la)=\left\{\begin{array}{ll}
\rang\textrm{Res}_\la ((2s-n)R(s)) & \textrm{ if }\la\not=\ndemi\\
\rang\textrm{Res}_\la R(s) & \textrm{ if }\la=\ndemi\end{array}\right.\]
It will also be useful later to know that $(n-2\la)R(\la)$ is holomorphic on $\{\Re(\la)=\ndemi\}$, this is 
a direct consequence of \cite[Lem. 4.1]{GZ}, the fact that $n^2/4$ is not an $L^2$ eigenvalue of $\Delta_g$
and \cite[Lem. 4.9]{PP}.
Moreover from \cite{MM}, $R(\la)\in \Psi_0^{-2,\la,\la}(\bar{X})+\Psi^{\la,\la}(\bar{X})$ thus its Schwartz kernel splits into
$R(\la;m,m')=R_1(\la;m,m')+R_2(\la;m,m')$ with (using coordinates (\ref{coordeclate}))
\[R_1(\la;m,m')=x^\la x'^\la K_\la(m,m'), \quad \beta^*R_2(\la;\rho,\rho',R,\omega,y')=\rho^\la{\rho'}^\la 
F_\la(\rho,\rho',R,\omega,y')\]
\begin{equation}\label{formres}
K_\la \in C^{\infty}(\bar{X}\x\bar{X}), \quad F_\la\in C^{\infty}(\bar{X}\x_0\bar{X}\setminus \textrm{diag})
\end{equation}
and $F_\la$ having a conormal singularity of order $-2$ at the interior diagonal $\textrm{diag}$.
We will call $\beta^*R_2$ the \emph{normal part} of the resolvent.
The boundary defining function $R$ and the coordinate $\omega$ actually depend on the choice of 
local coordinates we are using for the blow-up $\bar{X}\x_0\bar{X}$, we can actually show 
a parity regularity, in some sense, for the Taylor expansion of the normal part of the resolvent
at the front face $\mc{F}$, that is not depending on choice of coordinates.
To define regular parity, we use the special class of boundary defining function $x$ of $\bar{X}$, which 
induces $x:=\pi_L^*x$ and $x':=\pi_R^*x$ as geodesic boundary defining functions of $\bar{X}\x\bar{X}$.
The fibre $\mc{F}_p$ of the front face $\mc{F}$ (with $p=(y',y')\in M\x M$) 
being 
\begin{equation}\label{fp}
\mc{F}_p= \Big(\Big(N_p(\textrm{diag}_M,M\x M)\x(\rr^+\pl_x) \x (\rr^+\pl_{x'})\Big)
\setminus\{0\}\Big)/\{(w,t,u)\sim s(w',t',u'), s>0\}\end{equation}
there is an involution $\iota: (w,t,u)\to (-w,t,u)$ that passes to the quotient $\mc{F}_p$. 
Since $T_{y'}M$  is canonically isomorphic to $N_p(\textrm{diag}_M,M\x M)$ by  
$z\in T_{y'}M\to (z,-z)\in T_p(M\x M)$, $h_0(y')$ induces a metric on $N_p(\textrm{diag}_M,M\x M)$.
Then $\mc{F}_p$ is clearly identified with the quarter of sphere 
\[\mc{F}_p\simeq \{w+t\pl_x+u\pl_{x'}\in N_p(\textrm{diag}_M,M\x M)\x(\rr^+\pl_x) \x (\rr^+\pl_{x'}), t^2+u^2+|w|_{h_0(y')}^2=1\}\] 
this is actually trivial to check that the involution $\iota$ is just the symmetry $w\to -w$ in that model. 
In projective coordinates $(s=t/u,z=w/u)$, the interior of the face $\mc{F}_p$ is
a half-space diffeomophic to $\hh^{n+1}$ and the involution becomes $z\to-z$. 
Recall that $\beta^*x'$ is a smooth function that defines globally the interior of 
$\mc{F}$ in the sense that it vanishes and has non-degenerate differential there.

\begin{defi}\label{parity} 
One says that a function $L$ on $\mc{F}_p$ has even parity (resp. odd parity) if 
$L([-w,t,u])=L([w,t,u])$ (resp. $L([-w,t,u])=-L([w,t,u])$) in the model (\ref{fp})
where brackets denote equivalence classes.
A function $L\in \rho^{a}{\rho'}^bC^{\infty}(\bar{X}\x_0\bar{X})$ with $a,b\in\cc$ is said to have regular parity if 
for $x$ geodesic boundary defining function of $M$,  
then setting $x':=\beta^*\pi_R^*x$ with notations (\ref{beta})-(\ref{pirpil}), $L$ has the Taylor expansion at the interior of the front face $\mc{F}$
\[L=\sum_{i=0}^n{x'}^iL_i +O({x'}^{n+1}), \quad L_i\in \rho^a{\rho'}^{b-i} C^{\infty}(\mc{F}) \] 
with $L_{2i}$ having even parity and $L_{2i+1}$ odd parity on each $\mc{F}_p$. 
\end{defi} 
Note that one can define similarly regular parity for a conormal distribution to the 
interior diagonal $\textrm{diag}$, the $L_i$ will then be conormal distributions to 
$\pl\textrm{diag}:=\diag\cap\mc{F}$ on $\mc{F}$.
\begin{lem}\label{invariance} 
The regular parity is invariant with respect of choice of geodesic boundary
defining function $x$ of $\bar{X}$ as long as the metric $g$ on $X$ is even modulo $O(x^{n+1})$. 
Moreover, setting $x:=\beta^*\pi^*_Lx$ with notations (\ref{beta})-(\ref{pirpil}), the regular parity of $L\in \rho^{a}{\rho'}^bC^{\infty}(\bar{X}\x_0\bar{X})$ is equivalent to having the expansion 
\[L=\sum_{i=0}^nx^iL_i +O(x^{n+1}), \quad L_i\in {\rho}^{a-i}{\rho'}^b C^{\infty}(\mc{F}) \] 
with $L_{2i}$ having even parity and $L_{2i+1}$ odd parity on each $\mc{F}_p$.
\end{lem}
\textsl{Proof}: we consider a neighbourhood of $\mc{F}_p$
for some $p=(y_0,y_0)\in M\x M$. Let us take another geodesic boundary defining 
function $\hat{x}$ in $\bar{X}$, then from \cite{GR} one has that  
\[x=\hat{x}\sum_{2i\leq n-1}\hat{x}^{2i}f_{2i}+O(\hat{x}^{n+1}), \quad f_{2i}\in C^\infty(M).\] 
Thus we obtain, using notations $x'=\beta^*\pi^*_Rx$ and ${\hat{x}}'=\beta^*\pi_R^*\hat{x}$, 
\[x'={\hat{x}}'\sum_{2i\leq n-1}({\hat{x}}')^{2i}\beta^*\pi_R^*f_{2i}+O(({\hat{x}}')^{n+1})\]
and it clearly proves the first part of the Lemma since $\beta^*\pi_R^*f_{2i}$ are constant on the fibres of
$\mc{F}$ and only even powers of $\hat{x}'$ appear up to order $({\hat{x}}')^{n+1}$.
The second part can be checked through the change of projective coordinates near $\mc{F}_p$ 
\[(x',y',s=x/x',z=(y-y')/x')\to (x=x's,y=y'+x'z,t=1/s,z'=z/s),\]
where $x$ defines $\mc{F}$ and $(t,z')$ parametrize the fibers of $\mc{F}$ and $\iota$ is $z'\to -z'$; we have
\[\sum_{i=0}^n{x'}^iL_i(y',s,z)=\sum_{i=0}^n x^i t^iL_i(y-xz',1/t,z'/t)=
\sum_{i=0}^n x^iL'_{i}(y,t,z')+O(x^{n+1})\]
with $L'_j(y,t,z')=\sum_{i+|\alpha|=j}(\alpha!)^{-1}t^i(-z')^{\alpha}.\pl^\alpha_{y'}L_i(y,1/t,z'/t)$, this easily gives
$L_j'(y,t,-z')=(-1)^jL'_{j}(y,t,z')$ by asumption on $L_i$.
The equivalence holds by symmetry of the blow-up. 
\qed\\

\begin{remark}:  taking a geodesic boundary defining function $x$ and any coordinates $y$ near 
$y_0\in M$ we have coordinates $(x,y,x',y')$ with $(x,y)=\pi_L^*(x,y)$, $(x',y'):=\pi^*_R(x,y)$ 
near $p=(y_0,y_0)\in \textrm{diag}_M$ and induced coordinates in a neighbourhood of $\mc{F}_p$ 
\[R:=(x^2+{x'}^2+|y-y'|^2)^{\demi}, \quad \rho:=\frac{x}{R}, \quad \rho'=\frac{x'}{R}, 
\quad \omega=\frac{y-y'}{|y-y'|}, \quad y'\]
and since $x'=R\rho, s=\rho/\rho',z=\omega(1-\rho^2-{\rho'}^2)^{\demi}/\rho'$ are the projective coordinates on 
$\mc{F}$, it is clear that the regular parity property of a function 
$L\in \rho^{a}{\rho'}^bC^{\infty}(\bar{X}\x_0\bar{X})$ in this neighbourhood can be rephrased by
$L=\sum_{i=0}^nR^iL_i+O(R^{n+1})$ with  $L_i\in\rho^{a}{\rho'}^bC^{\infty}(\mc{F})$ and 
\[\iota^* L_i(\rho,\rho',\omega,y')=L_i(\rho,\rho',-\omega,y')=(-1)^iL_i(\rho,\rho',\omega,y').\]  
Clearly, the vector field $\beta^*(x\pl_x)$ acting on a function with regular parity
on $\bar{X}\x_0\bar{X}$ preserves its regular parity, whereas $\beta^*(x\pl_y)$ changes the parity.
\end{remark}

Since a function on $\mc{F}$ with odd parity restricts to $0$ at $\pl\textrm{diag}=\mc{F}\cap\textrm{diag}$ 
(the interior diagonal intersects $\mc{F}_p$ at a fixed point of the involution $\iota$
of $\mc{F}_p$), we get the straightforward 
\begin{lem}\label{paritediag}
Let $L\in\rho^a{\rho'}^b C^{\infty}(\bar{X}\x_0\bar{X})$ with $a,b\in\cc$ 
and suppose that $L$ has regular parity. Then the restriction of $F$ at the
diagonal $\textrm{diag}$ is a smooth function even modulo $O(x^{n+1})$ on $\bar{X}$ in the sense that
$m\in X\to \beta_*L(m,m)$ extends smoothly to $\bar{X}$ with an even expansion in powers of $x(m)$ modulo $O(x^{n+1})$ 
if $x$ is a geodesic boundary defining function.
\end{lem} 

A particularly interesting example of regular parity operator is the resolvent lifted kernel,
as we next prove in the  
\begin{prop}\label{transmres}
Let $(X,g)$ be an asymptotically hyperbolic metric  
which is even modulo $O(x^{n+1})$, then the normal part 
of the resolvent $R(\la)$ (see (\ref{formres})) has regular parity in the sense of Definition \ref{parity}, 
when $R(\la)$ is well-defined.
\end{prop} 
\textsl{Proof}:
to prove this result, one needs to return to the construction of the resolvent parametrix, in particular
the part involving the resolution on the front face, that is the normal operator.
The regular property is global in the sense that a fibre $\mc{F}_p$ of the front face 
is a manifold (with corners) but local in $\bar{X}\x_0\bar{X}$ since defined on 
each fibre, thus it suffices to work near an arbitrary fibre $\mc{F}_p$. 
Note that for $\Re(\la)>\ndemi$, (\ref{formres}) shows that the $n$ first asymptotic terms
of the lifted kernel of the resolvent at the front face $\mc{F}$ are given by those of its normal part.
If $p=(y_0,y_0)\in M\x M$, 
we have local coordinates $y$ near $y_0$ that induce projective coordinates 
\begin{equation}\label{proj}
(x',y',s=x/x',z=(y-y')/x')\end{equation}
near the interior $\{y'=y_0,x'=0)\}$ of $\mc{F}_p$ 
where again $x':=\beta^*\pi_R^*x$, $x:=\beta^*\pi_L^*x$, and similarly for $y,y'$.
Let us recall a couple of definitions and results about the normal operator, the 
reader could also read \cite{MM,JSB} or \cite[Sec. 2.3]{G1} for more details.
If $p=(y_0,y_0)\in\textrm{diag}_M$ with $y_0\in M$, the normal operator at $p$ 
of an operator $A\in\Psi_0^{m,k,l}(\bar{X})$ with Schwartz kernel $\kappa_A$ 
is the distribution $N_p(A):=\beta^*\kappa_A|_{\mc{F}_p}$ on $\mc{F}_p$ classically conormal 
to the point $\pl_p\textrm{diag}:=\textrm{diag}\cap \mc{F}_p$ 
of order $m$ and in $\rho^k{\rho'}^lC^\infty(\mc{F}_p\setminus\pl_p\textrm{diag})$, this distribution
can be interpreted as a left convolution distribution kernel on 
the Lie group $X_p:=\rr^+_s\x \rr^n_z$ with law $(s,z).(s',z'):=(ss',z+sz')$ and 
with the right invariant measure $|\det h_0(y_0)|^\demi s^{-1}dsdz$, 
thus as an operator on $X_p\simeq \hh^{n+1}$.  Note that the definition of $N_p(A)$ extends to the case 
$A\in R^\alpha\Psi_0^{m,k,l}(\bar{X})$ for $\Re(\alpha)>0$ by setting $N_p(A)=0$ (since the restriction 
on front face would vanish).
If $A\in\textrm{Diff}_0^d(\bar{X})\subset \Psi_0^{d,\infty,\infty}(\bar{X})$, its normal operator
is supported at $\pl_p\textrm{diag}$ and $N_p(\textrm{Id})=\delta_{\pl_p\textrm{diag}}$ is the Dirac mass there. 
If now $A\in\textrm{Diff}_0^d(\bar{X})$ and $B\in\Psi_0^{m,k,l}(\bar{X})$, then $AB\in\Psi_0^{m+d,k,l}(\bar{X})$  
and from \cite{MM} we have $N_p(AB)=N_p(A)N_p(B)$ in the sense of composition of 
the associated operators on $X_p$.
Thus if $R(\la)$ is the resolvent for $\Re(\la)>\ndemi$, we get
\begin{equation}\label{np}
N_p(\Delta_g-\la(n-\la))N_p(R(\la))=N_p(\textrm{Id})=\delta_{\pl_p\textrm{diag}}=\textrm{Id}_{X_p}\end{equation}
where again we identifies convolution kernel on $\mc{F}_p$ and operator on $X_p$.
From \cite{MM}, the normal operator $N_p(\Delta_g)$ is the Laplacian on the hyperbolic space $X_p\simeq \hh^{n+1}$,
which is easily seen from using projective coordinates (\ref{proj}) in (\ref{laplacian}) with $\beta^*(x\pl_x)=s\pl_s$, 
$\beta^*(x\pl_{y_i})=s\pl_{z_i}$
\begin{equation}\label{betadelta}
\beta^*\Delta_g=-(s\pl_s)^2+ns\pl_s+s^2\Delta_{h(x's,y'+x'z)}-\demi x's\tra_{h(x's,y'+x'z)}((\pl_xh)(x's,y'+x'z))s\pl_s,
\end{equation}
and freezing at $x'=0,y'=y_0$, that is $N_p(\Delta_g)=-(s\pl_s)^2+ns\pl_s+s^2\sum_{i,j}h^{ij}_0(y_0)\pl_{z_i}\pl_{z_j}$.
If one wants the solution of (\ref{np}) which is $L^2$ for $\Re(\la)>\ndemi$, $N_p(R(\la))$ is necessarily given by the resolvent of the Laplacian at energy $\la(n-\la)$ on $X_p$. In term of convolution kernel, thus of
distribution on $\mc{F}_p$, we have in projective coordinates (see \cite{MM,JSB} or \cite[App. B]{G1}) 
\begin{equation}\label{npr}
N_p(R(\la))(s,z)= G_\la\left(\frac{2s}{1+s^2+|z|_{h_0(y_0)}^2}\right)
\end{equation}
where $G_\la(x)=x^{\la}F_\la(x)$ for some $F_\la\in C^\infty([0,1))$ which can be expressed in term of
hypergeometric function. Note that the factor in $G_\la$ is nothing more than the inverse of the $\cosh$ of the 
hyperbolic distance of $(s,z)$ to $\pl_p\diag=\{s=1,z=0\}$ in $X_p$.
The result (\ref{npr}) holds for any $p=(y',y')$ if $y'$ is near $y_0$, with smooth dependence on $y'$ then. 
The first term $Q_{\la,0}(y',s,z):=\beta^*R(\la)|_{\mc{F}_{(y',y')}}$ of the (normal part of the) resolvent 
expansion at $\mc{F}$ clearly satisfies the regular parity property since it is even in $z$.
Locally $Q_{\la,0}$ can be considered as a kernel near $\mc{F}_p$, constant with respect to 
boundary defining function $x'$ of the interior of $\mc{F}$, it is the (local) lifted kernel
of an operator in $\Psi_0^{-2,\la,\la}(\bar{X})$.
One has locally near $\mc{F}_p$
\[\beta^*(\Delta_g-\la(n-\la))Q_{\la,0}=\delta_{\pl\textrm{diag}}+T_{\la,0}\]
where $T_{\la,0}$ is the (local) lifted kernel of an operator in $R\Psi_0^{0,\la,\la}(\bar{X})$ and actually 
in $R\Psi_0^{0,\la+1,\la}(\bar{X})$, see \cite{MM} for this fact.
Then we proceed by induction. Writing the lifted kernel of the normal part of the resolvent 
$\beta^*R_2(\la)$ near $\mc{F}_p$ as
\[\beta^*R_2(\la)=\sum_{k=0}^n{x'}^kQ_{\la,k}(y',s,z)+ O({x'}^{n+1}), \] 
with $Q_{\la,k}\in \rho^\la{\rho'}^{\la-k}C^\infty(\mc{F}\setminus \pl\textrm{diag})$, we assume
that for $k\leq j\leq n-1$ we have the regular parity
\begin{equation}\label{pariteq}
\iota^*Q_{\la,k}(y',s,z)=Q_{\la,k}(y',s,-z)=(-1)^kQ_{\la,k}(y',s,z)
\end{equation}
and we show that it also holds at order $j+1$. Let us denote by $P_\la$ the lift of $\Delta_g-\la(n-\la)$
under $\beta^*$, we then have 
\[P_\la\sum_{k=0}^j{x'}^kQ_{\la,k}=\delta_{\pl\textrm{diag}}+{x'}^{j+1}T_{\la,j}(x',y',s,z)\]
where $T_{\la,j}$ is the lifted kernel of an operator in $\Psi_0^{0,\la,\la-j}(\bar{X})$.
Since from (\ref{betadelta}), $P_\la$ commutes with $x'$, then it is clear that 
$P_\la|_{x'=0}Q_{\la,j+1}=T_{\la,j}|_{x'=0}$ and, like for the first step, this equation can be solved since $P_\la|_{x'=0}=N_{(y',y')}(\Delta_g-\la(n-\la))$ is the hyperbolic Laplacian. We solve it first 
near $\pl\textrm{diag}=\{s=1,z=0\}$ by symbolic parametrix up to continuous terms on $\mc{F}_p$ supported 
near $\pl\textrm{diag}$, this is standard and can be done by quantizing $\pl_s,\pl_z$, 
the key being ellipticity of $N_{(y',y')}(\Delta_g)$.   
The distribution $T_{\la,j}(0,y',s,z)$ is conormal of order $0$ at $\pl\textrm{diag}$
then by cutting it off near $\pl\diag$ by $\chi(2s/(1+s^2+|z|_{h_0(y')}^2))$ where 
$\chi\in C_0^{\infty}((1/2,2))$ equal $1$ near $1$, we have the oscillating integral
\[s^{-\ndemi}\chi\Big(\frac{2s}{1+s^2+|z|_{h_0(y')}^2}\Big)T_{\la,j}(0,y',s,z)=\int e^{i((s-1)\eta+z.\xi)}\sigma_{\la,j}(y';\eta,\xi)d\eta d\xi\]
for some $\sigma_{\la,j}$ which is a symbol of order $0$ in $(\eta,\xi)$. 
We will show later that
\begin{equation}\label{tlaj}
T_{\la,j}(0,y',s,-z)=(-1)^{j+1}T_{\la,j}(0,y',s,z),\end{equation} 
but let us assume it for now, then
$\sigma_{\la,j}(y';\eta,-\xi)=(-1)^{j+1}\sigma_{\la,j}(y';\eta,\xi)$ by inverse 
Fourier transform. We construct $\sigma^{(N)}_{\la,j+1}(y',s;\eta,\xi)$ by induction in $N$ by setting
$\sigma^{(0)}_{\la,j+1}(y',s;\eta,\xi):=\sigma_{\la,j}(y';\eta,\xi)$ 
and for $N\in\nn$
\[\sigma^{(N)}_{\la,j+1}(y',s;\eta,\xi):=L(s,\eta,\pl_s)\Big(\frac{\sigma^{(N-1)}_{\la,j+1}(y',s;\eta,\xi)}
{s^2\eta^2+s^2|\xi|^2_{h_0(y')}
+(\la-n/2)^2}\Big)\] 
where $L(s,\eta, \pl_s):=(s\pl_s)^2+2i\eta s(s\pl_s)+i\eta s$ and $\Re(\la)\gg n/2$. It is clear that 
$\sigma^{(N)}_{\la,j+1}$ is a symbol of order $-N$ in $(\eta,\xi)$, thus we can consider the oscillating integral 
\begin{equation}\label{qlaj}
Q^{(N)}_{\la,j+1}(y',s,z):=\chi\Big(\frac{2s}{1+s^2+|z|_{h_0(y')}^2}\Big)s^\ndemi \int e^{i(\eta (s-1)+\xi.z)}\sum_{k=0}^N\frac{\sigma^{(k)}_{\la,j+1}(y',s;\eta,\xi)}{s^2\eta^2+s^2|\xi|^2_{h_0(y')}
+(\la-n/2)^2}d\eta d\xi
\end{equation}
which is supported in a compact neighbourhood of $\pl\textrm{diag}$ in view of the cut-off. Since
\[s^{-\ndemi}N_{(y',y')}(\Delta_g)s^{\ndemi}=-(s\pl_s)^2-s^2\sum_{i,j}h^{ij}_0(y')\pl_{z_i}\pl_{z_j}+n^2/4,\] 
we get by construction that for $N$ large enough
\[(N_{(y',y')}(\Delta_g)-\la(n-\la))Q^{(N)}_{\la,j+1}(y',s,z)=
T_{\la,j}(0,y',s,z)-T^{(N)}_{\la,j}(y',s,z)\]
with $T^{(N)}_{\la,j}$ continuous on $\mc{F}$ near $\mc{F}_p$ 
(classically conormal of order $-N$ at $\pl\textrm{diag}$) and equal to 
$T_{\la,j}$ outside a neighbourhood of $\pl\textrm{diag}$. But we check that $Q^{(N)}_{\la,j+1}(y',s,-z)=(-1)^{j+1}Q^{(N)}_{\la,j+1}(y',s,z)$ 
using that $\sigma^{(k)}_{\la,j+1}$ has the parity of $\sigma_{\la,j}$ in $\xi$ (for any $k$) and 
the change of variable $\xi\to-\xi$ in (\ref{qlaj}).
Then we correct the error by setting
\[Q^{(\infty)}_{\la,j+1}(y',s,z):=\int_{\rr^+}\int_{\rr^n} G_\la\left(\frac{2u}{1+u+|v|_{h_0(y')}^2}\right)T^{(N)}_{\la,j}
\left(0,y',\frac{s}{u},z-\frac{s}{u}v\right)\frac{dudv}{u}|\det h_0(y')|^{\demi}\]
which is certainly convergent for $\Re(\la)\gg \ndemi$.
We clearly have $Q_{\la,j+1}=Q^{(N)}_{\la,j+1}+Q^{(\infty)}_{\la,j+1}$ and
by a change of variable $v\to -v$ in the last integral, the regular parity of $Q^{(\infty)}_{\la,j+1}$ holds 
if $T^{(N)}_{\la,j}(y',s,-z)=(-1)^{j+1}T^{(N)}_{\la,j}(y',s,z)$. This is actually straightforward 
in view of (\ref{tlaj}), the parity of $Q^{(N)}_{\la,j+1}$ and the fact that $N_{(y',y')}(\Delta_g)$
preserves parity in $z$ (it involves two differentiations in $z$).\\

It remains to prove (\ref{tlaj}), we first make an expansion of $P_\la$ 
at the front face $x'=0$. To that aim, we use the fact that 
$h(x)$ is even modulo $O(x^{n+1})$ and (\ref{laplacian}) to get
\[\Delta_g-\la(n-\la)=\sum_{2i\leq n}x^{2i}P_{2i,\la}(y;x\pl_x,x\pl_{y})+\sum_{2i+1\leq n}
x^{2i+1}P_{2i+1}(y;x\pl_{y})+O({x}^{n+1})\]
for some differential operators $P_{2i+1},P_{2i,\la}$ smooth in $y$, such that $Z\to P_{2i,\la}(y';Y,Z)$ 
is a polynomial even in $Z\in\rr^n$ (of degree $2$) and $Z\to P_{2i+1}(y';Z)$ is a polynomial
odd in $Z$ (of degree $1$). Since $y=y'+x'z$, $\beta^*(x\pl_x)=s\pl_s$, 
$\beta^*(x\pl_y)=s\pl_z$, the lift under $\beta$ gives (using multiindex $\alpha\in\nn^n$)
\[P_\la=\sum_{2i+|\alpha|\leq n}{x'}^{2i+|\alpha|}s^{2i}z^{\alpha}\pl^\alpha_{y'}
P_{2i,\la}(y';s\pl_s,s\pl_z)+\sum_{2i+1+|\alpha|\leq n}{x'}^{2i+1+|\alpha|}s^{2i+1}
z^{\alpha}\pl^\alpha_{y'}
P_{2i+1}(y';s\pl_z)\] 
modulo $O({x'}^{n+1})$. As a consequence,
\begin{eqnarray*}
T_{\la,j}(y',s,z)&=&\sum_{k+2i+|\alpha|=j+1}s^{2i}z^{\alpha}\pl^\alpha_{y'}P_{2i,\la}(y';s\pl_s,s\pl_z)
Q_{\la,k}(y',s,z)\\
 & &+\sum_{k+2i+|\alpha|=j}s^{2i+1}z^{\alpha}\pl^\alpha_{y'}P_{2i+1}(y';s\pl_z)
Q_{\la,k}(y',s,z).\end{eqnarray*} 
It is then easy to check the regular parity property from (\ref{pariteq}) combined with the parity properties of
$Z\to P_{2i,\la}(y';Y,Z)$, $P_{2i+1}(y';Z)$ and the fact that $(-1)^{k+|\alpha|}=(-1)^{k+|\alpha|+2i}$.

This achieves the proof by induction for $\Re(\la)\gg\ndemi$ and by 
meromorphic continuation in $\la$, the result extends as long as $\la$ is not a pole of $R(\la)$. 
\qed\\

\begin{remark}: Independently, Pierre Albin proved in \cite{A1} a similar result on
the renormalization of the resolvent on forms in order to analyze
the finite time $0$-Trace of the heat operator.\end{remark} 

\subsection{Eisenstein functions}
The Eisenstein functions are defined in this context by Joshi-Sa Barreto \cite{JSB} as boundary value
of the resolvent, or alternatively by 
Graham-Zworski \cite{GRZ} as solutions of a generalized Poisson problem. 
More precisely $E(\la)\in C^\infty(X\x M)
\cap C^{-\infty}(\bar{X}\x M)$ is the 
function, depending meromophically in $\la\in\cc\setminus -\nn$ if the metric
is even modulo $O(x^{n+1})$, defined by
\[E(\la):=2^{2\la-n}\frac{\Gamma(\la-\ndemi)}{\Gamma(\ndemi-\la)}\beta'_*(\beta^*({x'}^{-\la}R(\la))|_{\mc{B}})\]  
where we identified once again operator and Schwartz kernel. 
Since the kernel $R(\la;m,m')$ is symmetric in $(m,m')$, we have as a 
consequence of Proposition \ref{transmres} (and the remark following Lemma \ref{invariance})
and the definition of $E(\la)$ a regular parity for $E(\la)$ in the following sense 
\begin{cor}\label{transme}
If $x$ a geodesic boundary defining function of $\bar{X}$ and $y$ any coordinates 
near $y_0\in M$, then setting 
$(x,y):=\pi_L^*(x,y)$, $y':=\pi^*_Ry$, $R:= (x^2+|y-y'|)^{\demi}$,  
$\rho:=\frac{x}{R}$, $\omega:=\frac{y-y'}{R}$,  
\begin{equation}\label{ela}
E(\la;x,y,y')=\rho^{\la}R^{-\la}F_\la(\rho,R,\omega,y')+x^{\la}K_\la(x,y,y')\end{equation}
with $K_\la\in C^\infty(\bar{X}\x M)$ and $F_\la\in C^\infty(\bar{X}\x_0M)$ such that
$F_\la=\sum_{i=0}^nR^iF_{\la,i}+O(R^{n+1})$ with $F_{\la,i}\in 
C^{\infty}(\mc{F}')$ and
\begin{equation}\label{transmfi}
F_{\la,i}(\rho,-\omega,y')=(-1)^iF_{\la,i}(\rho,\omega,y').\end{equation}
\end{cor}

Let us recall the Green formula relating resolvent kernel and Eisenstein functions 
(see for instance \cite{G})
\begin{equation}\label{green}
R(\la;m,m')-R(n-\la;m,m')=
(n-2\la)\int_{\pl\bar{X}}E(\la;m,y')E(n-\la;m',y')\textrm{d}_{h_0}(y')\end{equation}
for $\la$ such that $R(\la),R(n-\la)$ are defined. 
Observe that this formula implies the smoothness of the kernel of $R(\la)-R(n-\la)$ in $X\x X$.

\subsection{Scattering operator}
The scattering operator in (\ref{scat}) is the same as usual, defined and studied in this context 
in \cite{JSB,GRZ} but renormalized with the $\Gamma$ factors 
to remove the infinite rank poles that it contains at $\ndemi+\nn$. 
From its definition in (\ref{scat}), we remark that $S(\la)$ depends on 
the geodesic boundary defining functions of $\bar{X}$ used to write the Poisson problem,
we get trivially that for another choice of geodesic boundary defining functions $\hat{x}$
with $\hat{x}=xe^{\omega}$ for $\omega\in C^{\infty}(\bar{X})$, the scattering operator obtained
in (\ref{scat}) is related to $S(\la)$ by the covariant rule
\begin{equation}\label{confrule}
\hat{S}(\la)=e^{-\la\omega_0}S(\la)e^{(n-\la)\omega_0}, \quad \omega_0=\omega|_{M}.
\end{equation}
It is also true \cite[Prop. 3.3]{GRZ} that $S(\la)$ is self-adjoint for $\la\in\rr$.
Since $x,\hat{x}$ are model defining functions correponding to respective conformal representative
$h_0=x^2g|_{TM}$ and $\hat{h}_0=e^{2\omega_0}h_0$ of the conformal infinity, $S(\la)$ is a 
conformally covariant operator on $(M,[h_0])$. It does depend on the whole manifold $X$
but it is proved by Graham-Zworski \cite{GRZ} that
\[P_k:=S\Big(\ndemi+k\Big)\in \textrm{Diff}^{2k}(M)\] 
is a differential operator of order $2k$ depending only on the $2k$ first derivatives 
$(\pl^j_xh(0))_{j\leq 2k}$ of the metric at the boundary, it is then local. 
In the case of an even 
Poincar\'e-Einstein manifold, $P_k$ is a natural conformal operator on $(M,[h_0])$, called $k$-th 
GJMS conformal Laplacian, defined by Graham-Jenne-Manson-Sparling \cite{GJMS} that generalizes Yamabe operator 
(which is $P_1$) and Paneitz operator (which is $P_2$).

It is proved in \cite{JSB,GRZ} that the kernel of $S(\la)$ is obtained by
\begin{equation}\label{kernels}
S(\la)=(2\la-n){\beta_\pl}_*({\beta'}^*(x^{-\la}E(\la))|_{\mc{T}'})=
2^{n+1-2\la}\frac{\Gamma(\la-\ndemi+1)}{\Gamma(\ndemi-\la)}{\beta_\pl}_*
(\beta^*(x^{-\la}{x'}^{-\la}R(\la))|_{\mc{B}\cap\mc{T}}).\end{equation}
Following \cite{G}, when the metric is even modulo $O(x^{n+1})$ (resp. modulo $O(x^{\infty})$), 
it is a meromorphic family of Fredholm operators on $M$ in 
$\la\in \cc\setminus{-\nn}$ (resp. $\la\in\cc$) if $n$ is odd.
In addition, $S(\la)\in\Psi^{2\la-n,0}(M)$ with principal symbol 
\begin{equation}\label{symbprinc}
\sigma_{\textrm{pr}}(S(\la))=|\xi|_{h_0}^{2\la-n}
\end{equation} 
where $h_0=x^2g|_{TM}$ is the representative of the conformal infinity associated to $x$.
On the essential spectrum of $\Delta_g$, represented by $\{\Re(\la)=\ndemi\}$, $S(\la)$ is 
a unitary operator satisfying
\begin{equation}\label{propscat}
S^{*}(\la)=S^{-1}(\la)=S(n-\la),\end{equation}
the second identity extends meromorphically in $\la\in\cc\setminus{-\nn}$ if $n$ is odd.  
A straightforward consequence of (\ref{transme}) and (\ref{kernels}) is 
\begin{prop}\label{sreg}
If $\la$ is not a singularity of $S(\la)$, 
then $S(\la)\in\Psi_{\rm reg}^{2\la-n,0}(M)$. If $\la$ is not a singularity
of $\pl_\la S(\la)S^{-1}(\la)$ then 
$\pl_\la S(\la)S^{-1}(\la)\in \Psi^{0,1}_{\rm odd}(M)$.
\end{prop}

Of course, using Proposition \ref{kvtrace}, this proves that $\TR(\pl_\la S(\la)S^{-1}(\la))$
is well defined.

\subsection{The main proof}
To obtain Krein's formula, it suffices actually to compute the $0$-Trace of the spectral measure
$d\Pi$ of $\Delta_g$, and by Stone's theorem, its Schwartz kernel satisfies
\[-2i\pi d\Pi(t^2;m,m')=R(\ndemi+it;m,m')-R(\ndemi-it;m,m'), \quad t\in(0,\infty)\] 
which is a smooth function in $X\x X$ by Green's formula (\ref{green}). First we check that
\begin{theo}\label{renorm}
Let $(X,g)$ be an asymptotically hyperbolic manifold of even dimension $n+1$  
whose metric is even modulo $O(x^{n+1})$. 
If $\la\notin \demi\zz$ is not a singularity of $R(\la),R(n-\la)$ and if $t\in(0,\infty)$, 
the values $\zerotr((n-2\la)(R(\la)-R(n-\la)))$ and $\zerotr(2td\Pi(t^2))$
are well defined and do not depend on $x$. 
\end{theo}
\textsl{Proof}: by smoothness of 
the Schwartz kernel of $(n-2\la)(R(\la)-R(n-\la))$ at the diagonal of $X\x X$ in view of (\ref{green}) 
we get that its lifted kernel by $\beta^*$ is in 
\[
(\rho\rho')^\la C^\infty(\bar{X}\x_0\bar{X})+(\rho\rho')^{n-\la}C^\infty(\bar{X}\x_0\bar{X})+
\beta^*((xx')^{\la}C^{\infty}(\bar{X}\x\bar{X}))+\beta^*((xx')^{n-\la}C^{\infty}(\bar{X}\x\bar{X}))
\]
and is holomorphic on the critical line $\Re(\la)=n/2$ since $(n-2\la)R(\la)$ is holomorphic 
on this line. 
Restricting the last two terms (their push-forward) at the diagonal of $X\x X$ gives 
the sum of a function in $x^{2\la}C^{\infty}(\bar{X})$ and one in $x^{2n-2\la}C^{\infty}(\bar{X})$
which renormalize independently of $x$ if $2\la\not\in\zz$ according to Subsection \ref{renormalized}.
The two first terms restricted to $\textrm{diag}$ are smooth, thus give a function $H\in C^\infty(\bar{X})$,
this function $H$ is the difference of the normal part of the resolvents $(n-2\la)R(\la)$ and 
$(n-2\la)R(n-\la)$, 
but Proposition \ref{transmres} shows that these normal parts have regular parity and 
Lemma \ref{paritediag} then gives us that $H$ has an even expansion in $x$ modulo $O(x^{n+1})$. 
We finally use result of Albin \cite[Th. 2.5]{A} to conclude independence with respect to $x$. 
The Stone formula gives the desired result for $2td\Pi(t^2)$ if $t>0$. 
\qed\\

Theorem \ref{renorm} also shows that $\zerotr(2td\Pi(t^2))$ extends analytically to $t\in\rr^*:=\rr\setminus\{0\}$ 
as an even function on $\rr^*$.

Next we compute this $0$-Trace in function of the scattering operator. The essential ingredients
are the Mass-Selberg relation (consequence of Green's formula (\ref{green})), the regular parity  
of the resolvent, Eisenstein function and scattering operator and the Proposition \ref{kvtrace}.

\begin{theo}\label{rs}
For $\la,n-\la\in\cc\setminus{\demi\zz}$ not resonance, we have 
\[\zerotr\Big((n-2\la)(R(\la)-R(n-\la))\Big)=\TR(\pl_\la S(\la)S^{-1}(\la)).\]
\end{theo} 
\textsl{Proof}: let us first change the notation and set $E(\la)$ for $(2\la-n)E(\la)$ so that
from (\ref{kernels}), its restricted kernel is the kernel of $S(\la)$. 
We fix a geodesic boundary defining function $x$ and use Mass-Selberg (following from (\ref{green})) 
identity like in Guillop\'e-Zworski \cite{GZ} or Patterson-Perry \cite{PP} to obtain
\[(2\la-n)\int_{x(m)>\eps}[R(\la;m,m')-R(n-\la;m,m')]_{m=m'}\textrm{d}_{g}(m)=(n-2\la)^{-1}(I_1(\la,\eps)-I_2(\la,\eps))\]
\[I_1(\la,x):=x^{-n}\int_{\pl\bar{X}}\int_{\pl\bar{X}}x\pl_x\pl_\la E(\la;x,y,y')E(n-\la;x,y,y')\textrm{d}_{h(x)}(y)\textrm{d}_{h_0}(y')\]
\[I_2(\la,x):=x^{-n}\int_{\pl\bar{X}}\int_{\pl\bar{X}}\pl_\la E(\la,x,y,y')x\pl_xE(n-\la;x,y,y')
\textrm{d}_{h(x)}(y)\textrm{d}_{h_0}(y').\]
Let us compute the finite part of $I_2(\la,x)$ as $x\to 0$, the one with $I_1(\la,x)$ being similar.
Let $(\chi_i)_i$ be a partition of unity of $M=\pl\bar{X}$ as in Proposition \ref{kvtrace} 
which we include in the integral $I_2$ to work in charts
\begin{eqnarray*}
I_2(\la,x)&=&\sum_{i,j}x^{-n}\int_{M\x M}\chi_i\otimes\chi_j\pl_\la E(\la)x\pl_xE(n-\la)
\textrm{d}_{h(x)}\otimes\textrm{d}_{h_0}\\
&=&\sum_{i,j}I^{ij}_{2}(\la,x).
\end{eqnarray*} 
This integral splits into two parts, one supported near the diagonal
of $M\x M$ (when $\chi_i\chi_j\not=0$) and the other supported far from the 
diagonal (when $\chi_i\chi_j=0$).
First suppose that $\chi_i\chi_j=0$, then 
\[\chi_i(y)\chi_j(y')E(\la;x,y,y')=x^{\la}\chi_i(y)\chi_j(y')K_\la(x,y,y')\]  
with $K_\la$ smooth and the integrand of $I^{ij}_{2}(\la,x)$ is 
\[\log(x)((n-\la)K_\la K_{n-\la}+xK_\la\pl_xK_{n-\la})+ (n-\la)\pl_\la K_\la K_{n-\la}+
x\pl_\la K_\la\pl_x K_{n-\la}.\]
Since the boundary value of $K_{\la}(x,y,y')|_{x=0}$ is $S(\la;y,y')$ and $K_\la$ smooth, 
one obtains directly that 
\begin{equation}\label{offdiag}
\textrm{FP}_{x= 0}I_2^{ij}(\la,x)=(n-\la)\int_{M\x M}\chi_i\otimes\chi_j\pl_\la S(\la)S(n-\la)
\textrm{d}_{h_0}\otimes\textrm{d}_{h_0}\end{equation}
and we are done for the off-diagonal part.\\

Now we have to deal with the part where $\chi_i\chi_j\not=0$, that is near the diagonal of $M\x M$.
Following (\ref{ela}) the kernel of $E(\la)$ is decomposed on $\textrm{supp}(\chi_i)\cap\textrm{supp}(\chi_j)$
into two parts
\[E(\la;x,y,y')=\rho^{\la}R^{-\la}F_\la(\rho,R,\omega,y')+x^{\la}K_\la(x,y,y')\]
with $K_\la\in C^\infty(\bar{X}\x M)$ and $F_\la\in C^\infty(\bar{X}\x_0M)$ in  
local coordinates on the blow-up $\bar{X}\x_0 M$ near the front face 
\begin{equation}\label{coord}
R=(x^2+r^2)^{\demi}, \quad r=|y-y'|,\quad \rho=\frac{x}{R}, \quad\omega=\frac{y-y'}{r}.
\end{equation}
An easy computation yields 
\begin{equation}\label{xdx}
x\pl_x=(1-\rho^2)\rho\pl_\rho+\rho^2R\pl_R,\end{equation}
therefore we get on $\textrm{supp}(\chi_i)\cap\textrm{supp}(\chi_j)$ 
\[x^{-n}\pl_\la E(\la)x\pl_xE(n-\la)= Q_{\textrm{sm}}(\la)+Q_{\sing}(\la)+
Q_\textrm{mix}(\la)\]
with 
\[Q_{\textrm{sm}}(\la):=\log(x)\Big((n-\la)K_\la K_{n-\la}+xK_\la\pl_xK_{n-\la}\Big)+ (n-\la)\pl_\la K_\la K_{n-\la}+
x\pl_\la K_\la\pl_x K_{n-\la}\]
the smooth term on the diagonal,  
\begin{eqnarray*}
Q_{\sing}(\la)&:=&R^{-2n}\Big((n-\la)(1-2\rho^2)(\pl_\la F_\la -\log(R^2)F_\la)F_{n-\la}+
(\pl_\la F_\la-\log(R^2)F_\la) x\pl_xF_{n-\la}\Big)\\
& &+ R^{-2n}\log(x)\Big((n-\la)(1-2\rho^2) F_\la F_{n-\la}+F_\la x\pl_xF_{n-\la}\Big)
\end{eqnarray*}
the singular term at the diagonal and 
\begin{eqnarray*}
Q_{\textrm{mix}}(\la)&:=&\log(x)\Big((n-\la)K_\la R^{-2n+2\la}F_{n-\la}+
K_\la x\pl_x(R^{-2n+2\la}F_{n-\la})+(n-\la)R^{-2\la}F_\la K_{n-\la}\\
& & +R^{-2\la}F_\la x\pl_xK_{n-\la}\Big)+ (n-\la)\pl_\la K_\la R^{-2n+2\la}F_{n-\la}+
\pl_\la K_\la x\pl_x (R^{-2n+2\la}F_{n-\la})\\
& & + (n-\la)\pl_\la(R^{-2\la}F_\la)K_{n-\la}+\pl_\la(R^{-2\la}F_\la)x\pl_xK_{n-\la}
\end{eqnarray*}
is the mixed term involving products of two types. 
Taking $(x^{-\la}E(\la))|_{\rho=0}$ in (\ref{ela}) the lifted kernel of $S(\la)$ can be decomposed under the form 
\[
\beta_\pl^*S(\la)=S_{\sing}(\la)+S_{\textrm{sm}}(\la)\]
\begin{equation}\label{ssing}
S_{\sing}(\la):= r^{-2\la}F_{\la}|_{\rho=0}\in r^{-2\la}C^\infty(M\x_0M), 
\quad S_{\textrm{sm}}(\la):=K(\la)|_{x=0}\in C^\infty(M\x M)\end{equation}
The finite part of the integral 
of $Q_{\textrm{sm}}(\la)$ in each chart is obtained directly like (\ref{offdiag}), that is, setting 
$\til{Q}_{\textrm{sm}}(\la):=\chi_i\otimes\chi_j Q_{\textrm{sm}}(\la)\textrm{d}_{h(x)}\otimes
\textrm{d}_{h_0}$,  
\begin{equation}\label{fpsm}
\textrm{FP}_{x=0}\int_{M\x M}\til{Q}_\textrm{sm}(\la)
=(n-\la)\int_{M\x M}\chi_i\otimes\chi_j(\pl_\la S_\textrm{sm}(\la) S_\textrm{sm}(n-\la))|_{x=0}
\textrm{d}_{h_0}\otimes\textrm{d}_{h_0}\end{equation}
Let us now consider the singular term $Q_{\sing}(\la)$ in the chart and on $\textrm{supp}(\chi_i\otimes\chi_j)$, it is clear that for $x<1$, $Q_\sing(\la)$ 
is supported in $R\in[0,A]$ for some $A>0$.
Let $\psi$ be a smooth cut-off function equal to $1$ in $[0,A]$ and $0$ in $[A+1,\infty)$. 
We use a Taylor expansion at $R=0$ to 
decompose $\til{Q}_{\sing}(\la):=(\chi_i\otimes \chi_j)Q_{\sing}(\la)\textrm{d}_{h(x)}\otimes\textrm{d}_{h_0}$ into
\[\til{Q}_\sing(\la)=\log(x)\Big(\sum_{i=0}^nR^{-2n+i}Q_{\log,i}(\la)+
R^{-n+1}Q_{\log,n+1}(\la)\Big)+\sum_{i=0}^nR^{-2n+i}Q_{i}(\la)+
R^{-n+1}Q_{n+1}(\la)\]
where, after identifying densities with functions via the local trivialisation $|dy\textrm{d}_{h_0}(y')|$, 
one has (recall $\mc{F}'=\{R=0\}$ is the front face of $\bar{X}\x_0 M$)
\begin{equation}\label{qlog}
Q_{\log,i}(\la;\rho,\omega,y')\in C^\infty(\mc{F}'), \quad i<n+1, \quad Q_{\log,n+1}(\la;\rho,R,\omega,y')\in 
C^\infty(\bar{X}\x_0 M)\end{equation}
\begin{equation}\label{qi}
Q_{i}(\la;\rho,R,\omega,y')=\log (R)Q^{\log}_i(\la)+Q^{\textrm{cst}}_i(\la)\in \log(R)C^\infty(\mc{F}')+C^\infty(\mc{F}'), \quad i<n+1
\end{equation} 
\begin{equation}\label{qn1}
Q_{n+1}(\la;\rho,R,\omega,y')\in \log(R)C^\infty(\bar{X}\x_0 M)+C^\infty(\bar{X}\x_0 M).\end{equation}
It is important to remark that the critical terms $Q_n,Q_{\log,n}$ are odd in the sense that
\begin{equation}
Q_{\log,n}(\la;\rho,-\omega,y')=-Q_{\log,n}(\la;\rho,\omega,y'), \quad Q_{n}(\la;\rho,R,-\omega,y')=-Q_{n}(\la;\rho,R,\omega,y'),
\end{equation}
this is a consequence of (\ref{transmfi}), $n$ is odd and the second remark following Lemma \ref{invariance}.
Observe that for $i<n$, a polar change 
of variable $(y,y')\to (u=|y-y'|/x, w=(y-y')/|y-y'|,y')$ gives, using (\ref{coord}),
\[\int \frac{Q_{\log,i}(\la)}{R^{2n-i}}=x^{-n+i}\int (1+u^2)^{-n+\frac{i}{2}}Q_{\log,i}\Big(\la;\frac{1}{(1+u^2)^{\demi}},
\frac{wu}{(1+u^2)^{\demi}},y'\Big)u^{n-1}dudw_{S^{n-1}}\textrm{d}_{h_0}(y')\]
\[\int \frac{Q_{i}(\la)}{R^{2n-i}}=x^{-n+i}\int (1+u^2)^{-n+\frac{i}{2}}Q^{\textrm{cst}}_{i}\Big(\la;\frac{1}{(1+u^2)^{\demi}},
\frac{wu}{(1+u^2)^{\demi}},y'\Big)u^{n-1}dudw_{S^{n-1}}\textrm{d}_{h_0}(y')\]
\[+x^{-n+i}\log x\int (1+u^2)^{-n+\frac{i}{2}}Q^{\log}_{i}\Big(\la;\frac{1}{(1+u^2)^{\demi}},
\frac{wu}{(1+u^2)^{\demi}},y'\Big)u^{n-1}dudw_{S^{n-1}}\textrm{d}_{h_0}(y')\]
\[+\demi x^{-n+i}\int (1+u^2)^{-n+\frac{i}{2}}\log(1+u^2)Q^{\log}_{i}\Big(\la;\frac{1}{(1+u^2)^{\demi}},
\frac{wu}{(1+u^2)^{\demi}},y'\Big)u^{n-1}dudw_{S^{n-1}}\textrm{d}_{h_0}(y')\]
all integrands are integrable and of the form $(C_1(\la)+C_2(\la)\log x)x^{-n+i}$ for constants $C_1(\la),C_2(\la)$ 
if $i<n$ in view of (\ref{qlog}) and (\ref{qi}). This also implies 
\begin{equation}\label{fpqlogi}
\textrm{FP}_{x=0}\log(x)\int \psi(R)R^{-2n+i}Q_{\log,i}(\la)=0
\end{equation}
\begin{equation}\label{fpqi}
\textrm{FP}_{x=0}\int \psi(R)R^{-2n+i}Q_{i}(\la)=\int (\psi(r)-1)Q_{i}(\la;0,r,w,y')r^{-n+i-1}drdw_{S^{n-1}}
\textrm{d}_{h_0}(y').
\end{equation}
The case $i=n$ is quite similar, we have for $x>0$
\[\int \frac{\psi(R)}{R^n}Q_{\log,n}(\la)=\int\frac{\psi((x^2+r^2)^\demi)}{(x^2+r^2)^{\ndemi}}
Q_{\log,n}\Big(\la;\frac{x}{(x^2+r^2)^{\demi}},
\frac{rw}{(x^2+r^2)^{\demi}},y'\Big)r^{n-1}drdw
\textrm{d}_{h_0}(y')\]
\[\int \frac{\psi(R)}{R^n}Q_{n}(\la)=\int\frac{\psi((x^2+r^2)^\demi)}{(x^2+r^2)^{\ndemi}}
Q_{n}\Big(\la;\frac{x}{(x^2+r^2)^{\demi}},(x^2+r^2)^{\demi},
\frac{rw}{(x^2+r^2)^{\demi}},y'\Big)r^{n-1}drdw
\textrm{d}_{h_0}(y')\]
and since we have proved in Subsection \ref{resolvent} that $Q_{\log,n}(\rho,-\omega,y')=-Q_{\log,n}(\rho,\omega,y')$ and the similar relation for $Q_n$, then
\begin{equation}\label{qn0}
\int \psi(R)R^{-n}Q_{\log,n}(\la)=\int \psi(R)R^{-n}Q_{n}(\la)=0, \quad \forall x>0
\end{equation}
thus the finite part is $0$. To deal with the remaining terms $Q_{\log,n+1},Q_{n+1}$, it suffices 
to use (\ref{qlog}) and Lebesgue theorem to see that
\[\log(x)\int \psi(R)R^{-n+1}Q_{\log,n+1}(\la)=\log(x)\int \psi(r)Q_{\log,n+1}(\la;0,r,w,y')drdw_{S^{n-1}}
\textrm{d}_{h_0}(y')+O(x^\demi)\]
whose finite part is $0$ and similarly from (\ref{qn1}) we get
\begin{equation}\label{fpreste}
\textrm{FP}_{x=0}\int \psi(R)R^{-n+1}Q_{n+1}(\la)=\int \psi(r)Q_{n+1}(\la;0,r,w,y')drdw_{S^{n-1}}\textrm{d}_{h_0}(y').
\end{equation}
The important fact deduced from (\ref{fpqlogi})-(\ref{fpreste}) is that
the finite part of the integral of $\til{Q}_{\sing}(\la)$ involves only the $(Q_i(\la)|_{\rho=0})_{i=0,\dots,n+1}$. 
But the $Q_i$ are defined by the polyhomogeneous expansion at $R=0$
\begin{equation}\label{dr2n}
AR^{-2n}\Big((n-\la)(1-2\rho^2)(\pl_\la F_\la -\log(R^2)F_\la)F_{n-\la}+
(\pl_\la F_\la-\log(R^2)F_\la) x\pl_xF_{n-\la}\Big)=\end{equation}
\[\sum_{i=0}^nR^{-2n+i}Q_i(\la)+Q_{n+1}(\la)R^{-n+1}=\sum_{i=0}^nR^{-2n+i}(\log (R)Q^{\log}_i(\la)+
Q_i^{\textrm{cst}}(\la))+Q_{n+1}(\la)R^{-n+1}\] 
where $A={\beta'}^*(\chi_i|\det h(x,y)|^\demi\otimes \chi_j)\in C^\infty(\bar{X}\x_0M)$.  Then the coefficients 
$(Q_i|_{\rho=0})_{i=0,\dots,n+1}$, in $\log(r)C^{\infty}(\mc{F}_\pl)+C^{\infty}(\mc{F}_\pl)$  
(recall $\mc{F}_\pl=\{r=0\}$ is the front face of $M\x_0 M$), are the coefficients 
in the expansion at $r=0$ of the restriction of (\ref{dr2n}) at $\rho=0$, which is 
\begin{equation}\label{nladrho}
(n-\la)A|_{\rho=0}r^{-2n}(\pl_\la F_\la|_{\rho=0} -\log(r^2)F_\la|_{\rho=0})F_{n-\la}|_{\rho=0}.\end{equation}
Returning to the definition (\ref{ssing}) of $S_{\sing}(\la)$, we deduce that the polyhomogeneous expansion
of
\[(n-\la)\beta_\pl^*(\chi_i|\det h_0(y)|^\demi\otimes\chi_j)\pl_\la S_{\sing}(\la)
S_{\sing}(n-\la)\]
at $r=0$ is (\ref{nladrho}).
Then combining (\ref{fpqlogi})-(\ref{fpreste}) we deduce
\begin{equation}\label{fpsing}
\textrm{FP}_{x=0}\int \til{Q}_{\sing}(\la)=
(n-\la)\int \Big\{(\psi(r)-1)\Big[\beta_\pl^*\big(\chi_i|\det h_0|^{\demi}\otimes\chi_j\big)
\pl_\la S_{\sing}(\la)S_{\sing}(n-\la)\Big]_{\sing}\end{equation} 
\[\quad\quad\quad+\psi(r)\Big[\beta_\pl^*\big(\chi_i|\det h_0|^{\demi}\otimes\chi_j\big)
\pl_\la S_{\sing}(\la)S_{\sing}(n-\la)\Big]_{\reg}\Big\}drdw_{S^{n-1}}\textrm{d}_{h_0}(y')\] 
with notations (\ref{usingreg}).
To conclude, it remains to study the mixed terms, and actually this can be done by analytic 
continuation in $\la\not\in\demi\zz$. Indeed for half of the terms (those with $R^{-2\la}$), 
one can take the limit in the integral (Lebesgue theorem) if $\Re(\la)>\ndemi$, this gives 
the limit as $x=0$ thus a finite part which admits a meromorphic continuation in $\la$, 
the other terms (those with $R^{2n-2\la}$) are dealt with similarly. As a consequence, we finally 
get, setting $\til{Q}_{\textrm{mix}}(\la):=\chi_i\otimes\chi_j Q_{\textrm{mix}}(\la)\textrm{d}_{h(x)}\otimes
\textrm{d}_{h_0}$,
\begin{equation}\label{fpmix}
\textrm{FP}_{x=0}\int \til{Q}_{\textrm{mix}}(\la)=(n-\la)\int 
\chi_i\otimes\chi_j[\pl_\la K_\la R^{2n-2\la}F_{n-\la}+\pl_\la(R^{-2\la}F_\la)K_{n-\la}]_{\rho=0}
\textrm{d}_{h_0}\otimes\textrm{d}_{h_0}
\end{equation}
\[\quad\quad\quad\quad=(n-\la)\int 
\chi_i\otimes\chi_j\Big(\pl_\la S_\textrm{sm}(\la) S_{\sing}(n-\la)+
\pl_\la S_{\sing}(\la)S_\textrm{sm}(n-\la)\Big)
\textrm{d}_{h_0}\otimes\textrm{d}_{h_0}\]
where we used (\ref{xdx}) again. Note that since the result holds for any cut-off function $\psi$ and is independent of 
the choice one can take the limit case $\psi=\indic_{[0,A]}$.
Using Proposition \ref{kvtrace} with $F:=\pl_\la S(\la)$ 
and $L:=S(n-\la)=S^{-1}(\la)$ which expresses the KV-Trace of $FL$
in term of lifted Schwartz kernels (note that we need $\trans S(n-\la)=S(n-\la)$, which
is a consequence of $R(n-\la)=\trans R(n-\la)$ and (\ref{kernels})) we see that
\[\textrm{FP}_{\eps\to 0} I_2(\la,\eps)=(n-\la)\TR(\pl_\la S(\la)S^{-1}(\la)).\] 

The part $I_1(x,\eps)$ is not really more complicated
and can be analyzed similarly. Using same notation as before, $Q_{\textrm{sm}}(\la),Q_{\sing}(\la),Q_{\textrm{mix}}(\la)$
become now
\[Q_{\textrm{sm}}(\la):=\Big(K_\la+\la \pl_\la K_\la+x\pl_x\pl_\la K_{\la}+\log(x)x \pl_xK_\la\Big)K_{n-\la}\]  
\begin{eqnarray*}
Q_{\sing}(\la)&:=&R^{-2n}\Big((1-2\rho^2)((1-\la\log(R^2))F_\la+\la\pl_\la F_\la-\log(R^2)x\pl_xF_\la\Big)F_{n-\la}\\
& &+R^{-2n}\log x\Big(\la(1-2\rho^2) F_\la +x\pl_x F_\la\Big)F_{n-\la}
\end{eqnarray*} 
\begin{eqnarray*}
Q_{\textrm{mix}}(\la)&:=&R^{-2n+2\la}\Big(K_\la +\la\pl_\la K_\la+x\pl_x K_\la +\log (x)x\pl_xK_\la\Big)F_{n-\la}\\
&&+R^{2\la}\Big( (1-2\rho^2)F_\la+\la(1-\rho^2)(\pl_\la F_\la -\log(R^2)F_\la)-\log(R^2)x\pl_xF_\la \Big)K_{n-\la}\\
&&+R^{2\la}\log x\Big(\la(1-2\rho^2)F_\la+x\pl_x F_\la \Big)K_{n-\la}.
\end{eqnarray*}
Applying exactly the same method, it is straightforward to see that $\la\TR(\pl_\la S(\la)S(n-\la))$ contributes to 
the finite part of $I_1(\la,\eps)$ as $\eps\to 0$ but it turns out that another term appears, 
this is the KV-Trace on $M$ of the composition of operators with respective lifted Schwartz kernel 
\[(r^{-2\la}F_\la+K_{\la})|_{\rho=0}, \quad (r^{-2n+2\la}F_{n-\la}+K_{n-\la})|_{\rho=0},\]
these operators are clearly $S(\la)$ and $S(n-\la)$ thus using $S(\la)S(n-\la)=\textrm{Id}$ and the fact that
the KV-density of Id is $0$ (the full symbol is constant), we see that this term does not contribute to the finite part
of $I_1(\la,\eps)$ as $\eps\to 0$. 
The proof of the Theorem is complete.
\qed\\
 
Since $\pl_\la S(\la)S^{-1}(\la)$ is meromorphic in $\la\in\cc\setminus (-\nn\cup n+\nn)$ (resp. $\cc$) 
and holomorphic on $\{\Re(\la)=\ndemi\}$ with values in $\Psi^{0,1}_\odd(M)$ according to Proposition \ref{sreg}
when $g$ is even modulo $O(x^{n+1})$ (resp. even modulo $O(x^\infty)$), then the formula of Lemma \ref{kvkernel} 
clearly implies that 
$\TR(\pl_\la S(\la)S^{-1}(\la))$ is meromorphic in $\cc\setminus (-\nn\cup n+\nn)$ (resp. $\cc$) 
and holomorphic on the line $\{\Re(\la)=\ndemi\}$.
Thus $\zerotr(2td\Pi(t^2))$ extends analytically to $t\in\rr$ and 
meromorphically to $\cc\setminus (i(n/2+\nn)\cup i(-n/2-\nn))$ (resp. to $\cc$)
if $g$ is even modulo $O(x^{n+1})$ (resp. even modulo $O(x^\infty)$).
We can then define the \emph{generalized Krein function} by
\begin{equation}\label{kreinfct}
\xi(t):=\int_{0}^t \zerotr (2ud\Pi(u^2))du
\end{equation}
analytic in $t\in\rr$, and odd on $\rr$. 
Formally, $\xi(t)$ is the $0$-Trace of $\Pi(t^2)$.
We have the identity
\begin{equation}\label{plxi}
\TR(\pl_\la S(\la)S^{-1}(\la))=-2\pi\pl_z\xi(z)|_{z=i(\ndemi-\la)}.
\end{equation}
Let us now study the singularities of $\TR(\pl_\la S(\la)S^{-1}(\la))$.

\begin{prop}\label{residus}
The function $s(\la):=\TR(\pl_\la S(\la)S^{-1}(\la))$ has only first order poles
in the set $\cc\setminus (-\nn\cup n+\nn)$ with residues 
\[-{\rm Res}_{\la_0}s(\la)=m(\la_0)-m(n-\la_0)+\indic_{\ndemi-\nn}(\la_0)\dim \ker S(n-\la_0)-
\indic_{\ndemi+\nn}(\la_0)\dim\ker S(\la_0)\] 
the result holds in $\cc$ if $g$ is even modulo $O(x^{\infty})$.
\end{prop}
\textsl{Proof}: we define as in \cite{G2} the operator
\begin{equation}\label{defstil}
\til{S}(\la):=P^{-\frac{\la}{2}+\frac{n}{4}}S(\la)P^{-\frac{\la}{2}+\frac{n}{4}}, \quad P:=\Delta_{h_0}+1
\end{equation}
where $h_0=x^2g|_{TM}$ if $S(\la)$ has been defined using geodesic boundary defining function $x$.
Then using Lemma \ref{power} and (\ref{symbprinc}) we have $\til{S}(\la)\in\Psi_\odd^{0}(M)$
and $\til{S}(\la)-1\in\Psi^{-1}(M)$ is a family of compact operators.
We have from Gohberg-Sigal theory (see \cite{GS} or \cite[Eq. 2.2]{G2}) the factorization 
near $\la_0$
\begin{equation}\label{factorization}
\til{S}(\la)=U_1(\la)\Big(P_0+\sum_{l=1}^m(\la-\la_0)^{k_l}P_l\Big)U_2(\la)
\end{equation}
where $U_1,U_2$ are holomorphically invertible bounded operators on $L^2(M)$, $P_i$
some orthogonal projectors with $\rang P_l=1$ if $l>0$, $\sum_{l=0}^mP_l=1$, $P_iP_j=\delta_{ij}P_i=\delta_{ij}P_j$, 
and $k_l\in\zz\setminus \{0\}$.
The inverse has to be 
\[\til{S}^{-1}(\la)=U_2(\la)^{-1}\Big(P_0+\sum_{l=1}^m(\la-\la_0)^{-k_l}P_l\Big)U_1(\la)^{-1}.\]
Moreover, combining \cite[Eq 2.4]{G2}\footnote{there is a sign typo for the definition
of $N_{\la_0}(M^{-1}(\la))$!} with \cite[Th. 1.1]{G2}, completed by \cite[Prop. 2.1]{GN} for points of $\sigma_{\textrm{pp}}(\Delta_g)$), we get
\begin{equation}\label{sumkl}
-\sum_{l=1}^mk_l=m(\la_0)-m(n-\la_0)+\indic_{\ndemi-\nn}(\la_0)\dim \ker S(n-\la_0)-
\indic_{\ndemi+\nn}(\la_0)\dim\ker S(\la_0).\end{equation}
and it remains to show that ${\rm Res}_{\la_{0}}s(\la)=\sum_{l=1}^mk_l$.
Let us first compute 
\begin{equation}\label{tilvspastil}
\pl_\la S(\la)S^{-1}(\la)=\demi\log(P)+P^{\frac{\la}{2}-\frac{n}{4}}\pl_\la\til{S}(\la)\til{S}^{-1}(\la)
P^{-\frac{\la}{2}+\frac{n}{4}}+\demi
P^{\frac{\la}{2}-\frac{n}{4}}\til{S}(\la)\log(P)\til{S}^{-1}(\la)P^{-\frac{\la}{2}+\frac{n}{4}}
\end{equation}
Using Lemma \ref{power} with $P$ (which has regular parity since differential of order $2$), 
all the terms that appear have regular parity, this means from 
(\ref{traceab}) that we can use cyclicity of the KV-Trace and deduce
\begin{equation}\label{tilpastil}
\TR(\pl_\la S(\la)S^{-1}(\la))=\TR(\log(P))+\TR(\pl_\la\til{S}(\la)\til{S}^{-1}(\la)).
\end{equation}
Now from (\ref{factorization}) we get
\[\pl_\la\til{S}(\la)\til{S}^{-1}(\la)=(\la-\la_0)^{-1}\sum_{l=1}^m k_lU_1(\la)P_lU_1^{-1}(\la)+
\pl_\la U_1(\la)U_1^{-1}(\la) \]
\[\quad\quad\quad\quad\quad+\sum_{l,j=0}^m(\la-\la_0)^{k_l-k_j}U_1(\la)P_l\pl_\la U_2(\la)U_2(\la)^{-1}P_jU_1^{-1}(\la)\]
where we have set $k_0=0$. The polar part has finite rank (the term $l=j=0$ is holomorphic at $\la_0$)
thus of trace class, this implies that the KV-Trace (viewed as $\TR(A):=\tra(AP^s)|_{s=0}$ and
not necessarily on $\psdo$'s) 
\begin{equation}\label{TRhol}
\TR\Big(U_1(\la)P_0\pl_\la U_2(\la)U_2(\la)^{-1}P_0U_1^{-1}(\la)+\pl_\la U_1(\la)U_1^{-1}(\la)\Big)
\end{equation}
is well-defined and analytic, at least for $\la\not=\la_0$.
Using cyclicity of the trace $\tra(AB)=\tra(BA)$ if $A$
is bounded and $B$ trace class, we get that the trace, or KV-Trace, of the polar part is the sum of the following 
two terms
\begin{equation}\label{trpolar1}
\TR\Big((\la-\la_0)^{-1}\sum_{l=1}^m k_lU_1(\la)P_lU_1^{-1}(\la)\Big)=(\la-\la_0)^{-1}\sum_{l=1}^m k_l,
\end{equation}
\begin{equation}\label{trpolar2}
\TR\Big(\sum_{\substack{l+j>0\\ k_l-k_j<0}}^m(\la-\la_0)^{k_l-k_j}U_1(\la)P_l\pl_\la U_2(\la)U_2(\la)^{-1}P_jU_1^{-1}(\la)\Big)
=0,\end{equation} 
the last identity coming also from $P_lP_j=\delta_{lj}P_j$.
Thus the term (\ref{TRhol}) is meromorphic, it remains to prove that its polar part at $\la_0$ 
is $0$. This can be done easily by setting $Z(\la)$ the polar part (with finite rank) of 
$\pl_\la\til{S}(\la)\til{S}^{-1}(\la)$ and checking that for $\Re(s)\ll 0$ and $i\in\nn_0$
\[\int_{C(\la_0)}(\la-\la_0)^i\tra((\pl_\la\til{S}(\la)\til{S}^{-1}(\la)-Z(\la))P^s)d\la=
\tra\int_{C(\la_0)}(\la-\la_0)^i(\pl_\la\til{S}(\la)\til{S}^{-1}(\la)-Z(\la))P^sd\la\]
is $0$ by holomorphy of $\pl_\la\til{S}(\la)\til{S}^{-1}(\la)-Z(\la)$, here $C(\la_0)$
is a little circle around $\la_0$. As a conclusion (\ref{TRhol}) is holomorphic at $\la_0$ and combining
this fact with (\ref{tilpastil}), (\ref{trpolar1}), (\ref{trpolar2}), (\ref{sumkl}),
the proof is achieved. 
\qed\\

For $g$ even Poincar\'e-Einstein manifold, if $n^2/4-k^2\notin \sigma_{\textrm{pp}}(\Delta_g)$ then $\dim\ker S(n/2+k)$ 
are conformal invariants depending only on the conformal infinity since $S(n/2+k)=P_k$ is the $k$-th GJMS conformal 
Laplacian of $(M,[h_0])$.\\

A particular case of interest is when the curvature is constant, that is essentially
when $X=\Gamma\backslash\hh^{n+1}$ is a quotient of the hyperbolic space by
a convex co-compact group. There is in this case a dynamical zeta function introduced by Selberg
\begin{equation}\label{zeta}
Z(\la)=\exp\left(-\sum_{\gamma}\sum_{m=1}^{\infty}\frac{1}{m}\frac{e^{-\la ml(\gamma)}}{G_\gamma(m)}\right)\end{equation}
where $\gamma$ runs over the primitive closed geodesics of $X$,
$l(\gamma)$ is the length of $\gamma$ and $G_\gamma(m):=e^{-\ndemi ml(\gamma)}|\det(1-P_\gamma^m)|^{\demi}$ 
if $P_\gamma$ is the Poincar\'e linear map associated to the primitive periodic orbit $\gamma$ of the geodesic
flow on the unit tangent bundle. It is well-known that the infinite sum converges for $\Re(\la)>\delta$ if $\delta$ is the dimension the limit set of the group $\Gamma$, or equivalently the exponent of convergence of the Poincar\'e series 
of the group.
 
We then clearly deduce for this case, using Patterson formula (see e.g. the proof of Theorem 1.1 of \cite{GN}),
\begin{cor}\label{convexcc} 
If $(X,g)$ is a convex co-compact hyperbolic quotient $X=\Gamma\backslash \hh^{n+1}$ with $n$ odd, we then have 
\begin{eqnarray*}
\pl_z\xi(z)&=&\frac{i}{2\pi}\TR\Big(\pl_zS\Big(\ndemi+iz\Big)S^{-1}\Big(\ndemi+iz\Big)\Big)\\
&=&\frac{1}{2\pi}\Big(\frac{Z'(\ndemi+iz)}{Z(\ndemi+iz)}+\frac{Z'(\ndemi-iz)}{Z(\ndemi-iz)}+\pi^{-\ndemi}\frac{\Gamma(\ndemi)}{\Gamma(n)}\frac{\Gamma(\ndemi+iz)\Gamma(\ndemi-iz)}{\Gamma(iz)\Gamma(-iz)}
\zerov(X)\Big)\end{eqnarray*}
where $Z(\la)$ the Selberg zeta function of $\Gamma$, $\zerov(X)$ is the $0$-volume of $X$ -i.e. the $0$-integral of $1$- which is equal by \cite[Appendix]{PP} to $\zerov(X)=(-1)^{\frac{n+1}{2}}\pi^{\ndemi+1}\chi(\bar{X})/\Gamma(\ndemi+1)$ with $\chi(\bar{X})$ the Euler characteristic of $\bar{X}$. 
\end{cor}

This corollary together with Proposition \ref{residus} (and the computation 
by \cite[Appendix]{PP} of the $0$-volume of $X=\Gamma\backslash\hh^{n+1}$ in term of the 
Euler characteristic) gives another proof of Theorem 1.5 of \cite{PP} when $n$ is odd.
\begin{cor}\label{convexccc2}
The Selberg zeta function $Z(\la)$ for a convex co-compact quotient $X=\Gamma\backslash\hh^{n+1}$ 
in odd dimension has for divisors the resonances, some conformal invariants at $n/2-\nn$ 
and some topological invariants at $-\nn_0$ as written in \cite[Th. 1.5]{PP}
\end{cor} 
 
\section{Applications to determinants}

\subsection{Regularized determinant of $S(\la)$}
In this part we define the determinant of $S(\la)$ for $n$ odd and for even
modulo $O(x^{\infty})$ asymptotically hyperbolic metrics. 
Since $S(\la)$ is an elliptic self-adjoint classical $\psdo$ of order $2\la-n$ for $\la\in(\ndemi,+\infty)$,
one can define its regularized determinant using the method of Kontsevich-Vishik \cite[Sec. 2]{KV}, itself inspired from 
the zeta regularized determinant of Laplacians introduced by Ray-Singer \cite{RS}.
Let us first take $\la\not\in\ndemi+\nn$. If $S(\la)$ is invertible (which is the case for almost 
any $\la\in(\ndemi,+\infty)$) the logarithm of $S(\la)$ and the complex power $S(\la)^s$ for $s\in\cc$ can be constructed using a spectral cut of $S(\la)$, here for instance any 
$L_\theta:=\{re^{i\theta}, r\in(0,\infty))\}$ works as long as $\theta\not=0(\pi)$.
It then allows to define the spectral zeta function of $S(\la)$ by 
\[\zeta_{\la}(s):=\TR(S(\la)^{-s})\] 
where the trace is the KV-Trace. In view of \cite[Prop. 3.4]{KV} this function is meromorphic in $s\in\cc$ with  
at most simple poles at $(2\la-n)^{-1}(-n+\nn)$, the residue at a pole $s$ is given by the 
Wodzicki residue of $S(\la)^{-s}$. It is actually holomorphic at $s=0$ 
since Wodzicki residue of $\textrm{Id}$ vanishes and one can set
\begin{equation}\label{defdet}
\det S(\la):=e^{-\pl_s\zeta_\la(0)}.
\end{equation}
Note now that for $2\la-n=2k\in 2\nn$, $S(\la)$ is a differential operator which is in the 
odd class $\Psi^{2k}_{\odd}(M)$, then these definitions go through as well.
\begin{prop}\label{dets1}
The function $\la\to \det S(\la)$ is analytic for $\la\in (n/2,\infty)$ and extends meromorphically 
to $\la\in\cc$. It can be expressed under the form
\[\det S(\la)= \exp(\TR \log S(\la))\]
and its derivative is $\pl_\la\det S(\la)=\det S(\la) \TR(\pl_\la S(\la)S^{-1}(\la))$.   
\end{prop}
\textsl{Proof}: Let $\la_0>n/2$ such that $S(\la_0)$ is invertible and choose $\la>0$ in a small neighbourhood of $\la_0$, 
so that $S(\la)$ is invertible.
We now apply Corollary 2.4 of \cite{PS} with $Q:=S(\la)$ and $A:={\rm Id}$ in their notation. 
The local residue of $\log Q$ is defined by
 \[{\rm res}_y(\log Q)=\int_{S^{*}(\pl\bar{X})} \sigma_{-n,0}(\log Q)(y,\xi)d\xi, \quad y\in\pl\bar{X}\]
where $S^{*}(\pl\bar{X})$ is the cosphere bundle of $\pl\bar{X}$, 
$\sigma(\log Q)(y,\xi)$ is the local symbol of $\log Q$ with polyhomogeneous expansion
(i.e. $\sigma_{-j,k}(\log Q)(y,\xi)$ below is homogeneous of degree $-j$ in $\xi$)
\[\sigma(\log Q)(y,\xi)\sim\sum_{j=0}^\infty\sum_{k=0}^1\sigma_{-j,k}(\log Q)(y,\xi)(\log|\xi|)^{k}\]
This local residue vanishes 
since $\log S(\la)\in\Psi^{0,1}_{{\rm odd}}(\pl\bar{X})$ by Lemma \ref{power} and the integral becomes the integral of an odd
density on the sphere, so Corollary 2.4 of \cite{PS} reads
\[-\pl_s\zeta_{\la}(s)|_{s=0}=\TR(\log S(\la)).\] 
By construction, $\log S(\la)\in\Psi^{0,1}_{\odd}(\pl\bar{X})$ depends analytically on $\la$ in the sense of \cite[Def. 1.9]{PS} for $\la$ in a small real neighborhood of 
$\la_0$ (so that $S(\la)$ is self-adjoint, invertible, with positive principal symbol), then Lemma \ref{power} shows that $\pl_s\zeta_\la(0)$ is analytic in $\la$ near $\la_0$. Differentiating with respect to 
$\la$ gives 
\begin{equation}\label{dlogdetsla}
\pl_\la\log \det S(\la)=\TR(\pl_\la S(\la)S(\la)^{-1}).
\end{equation}
The right hand side extends to $\cc$ meromorphically by Proposition \ref{residus}, with 
first order poles only and whose residues are integers, this implies that
one can integrate (\ref{dlogdetsla}) and take the exponential, to get a meromorphic continuation of $\det S(\la)$
\[\det S(\la)=(\det S(\la_0))\exp\Big(\int_{\la_0}^\la\TR(\pl_zS(z)S^{-1}(z))dz\Big)\]
where $\la_0\in(\ndemi,\infty)$ is chosen such that $S(\la_0),S(\la_0)^{-1}$ exists.\qed\\

From (\ref{plxi}) one also obtains for any $\la_0=\ndemi+iz_0>\ndemi$ which is not singularity of $S(\la)$
\[\det S\Big(\ndemi+iz\Big)=e^{-2\pi i(\xi(z)-\xi(z_0))}\det S\Big(\ndemi+iz_0\Big).\]
It is also quite easy to see that $\det S(\la)$ is a conformal invariant of the conformal infinity 
$(M,[h_0])$. Indeed if $\hat{h}_0=e^{2\omega_0}h_0$ is another conformal representative of $[h_0]$,
one has another geodesic boundary defining function $\hat{x}$ in $\bar{X}$ which induces a scattering operator 
$\hat{S}(\la)$ related to $S(\la)$ by (\ref{confrule}), but the same works by taking 
$h^t_0:=e^{2t\omega_0}h_0$ for $t\in[0,1]$, that is we have a smooth one parameter 
family of associated scattering operators 
\[S_t(\la)=e^{-t\la\omega_0}S(\la)e^{(n-\la)t\omega_0}, \quad S_1(\la)=\hat{S}(\la).\]
all invertible if $S(\la)$ is. 
For $\la>\ndemi$ fixed such that $S(\la)$ is invertible, let us compute the logarithmic derivative with respect to $t$ using
method previsously detailed (cyclicity of $\TR$ for operator with regular parity, 
recalling that a multiplication operator has regular parity)
\[\pl_t\log\det S_t(\la)=-2\TR(\omega_0)=0\]  
the fact that it vanishes is an easy consequence of Lemma \ref{kvkernel} with 
the fact that $\omega_0$ has a symbol $a(y,\xi)=\omega_0(y)$ constant in $\xi$, thus its Schwartz kernel
is supported on the diagonal. We conclude that $\pl_t\det S_t(\la)=0$ since $\det S_t(\la)\not=0$
by assumption (it is defined as an exponential), then $\det \hat{S}(\la)=\det S(\la)$. 

Combining this discussion and Proposition \ref{residus} we obtain
\begin{prop}\label{mainth}
On an asymptotically hyperbolic manifold even modulo $O(x^{\infty})$ in even dimension,
the function $\det S(\la)$ has a meromorphic continuation to $\cc$
with divisor at $\la_0\in\cc$ given by
\[-m(\la_0)+m(n-\la_0)-\indic_{\ndemi-\nn}(\la_0)\dim \ker S(n-\la_0)+
\indic_{\ndemi+\nn}(\la_0)\dim\ker S(\la_0).\]
Moreover $\det S(\la)$ is a conformal invariant of the conformal infinity and satisfies 
\[\det S\Big(\ndemi+iz\Big)=Ce^{-2\pi i\xi(z)}\]
where $\xi$ is the extension in $\cc$, defined modulo $\zz$, 
of Krein's function introduced in (\ref{kreinfct}) and $C\in\cc$ constant.
\end{prop}

Remark that this allows to consider $\xi$ as a scattering phase, 
i.e a phase of the scattering operator.
We now want to consider $\det\til{S}(\la)$ where $\til{S}(\la)$ was defined in (\ref{defstil}), in order to compute
the contant $C$ of last Theorem.
Before defining this determinant, recall that
the construction of $\det A$ for an elliptic self-adjoint operator with positive principal symbol 
$A\in\Psi^{a}(M)$ with $a\in\rr_+$ can be defined as we did for $S(\la)$ (see again
\cite{KV}) by
\[\det A:=e^{-\pl_s\textrm{TR}(A^{-s})|_{s=0}}.\]
If the order $a$ is negative it is still possible to define the complex powers, 
thus the spectral zeta function, by $A^{-s}:=(A^{-1})^s$ where now $A^{-1}$ has order $-a>0$.
Since $\til{S}(\la)\in\Psi_\odd^{0}(M)$ with principal symbol equal to $1$, 
one can use the trick of Kontsevich-Vishik \cite[Cor. 4.1]{KV} to define its determinant, 
namely for any $\la$ such that $\til{S}(\la)$ is invertible, set
\begin{equation}\label{dettils}
\det \til{S}(\la):=\frac{\det(\til{S}(\la)P)}{\det P}\end{equation}
here $P:=(1+\Delta_{h_0})$ is clearly odd since it is differential and 
the dimension of $M$ is odd, thus $\til{S}(\la)P$ too. It is proved in \cite[Cor. 4.1]{KV} 
that the obtained determinant does not depend on $P$, is also equal to
\begin{equation}\label{propdet}
\det \til{S}(\la)=\frac{\det(P\til{S}(\la))}{\det P},\quad \det \til{S}(\la)=\frac{\det(\til{S}(\la)P^l)}{\det P^l}
\end{equation} 
for any power $l\in\nn$.
We wish to relate $\det\til{S}(\la)$ to $\det S(\la)$, to that aim we need a multiplication 
result for determinant due to Kontsevich-Vishik \cite[Th. 4.1]{KV}, or actually 
a generalization. 
\begin{lem}\label{multdet}
Let $\alpha$, $\beta$ be two non-zero real numbers with $\alpha+\beta\not=0$. 
If $A\in\Psi_\reg^{\alpha,0}(M), B\in\Psi^{\beta,0}_\reg(M)$ are elliptic self-adjoint invertible
with positive principal symbol, then 
\[\det(AB)=\det(A)\det(B).\]
\end{lem}
\textsl{Proof}: we follow the proof of Theorem 4.1 in \cite{KV}. First set for $t\in[0,1]$ 
\[\eta:= AB^{-\frac{\alpha}{\beta}}, \quad A_t:=\eta^tB^{\frac{\alpha}{\beta}}\]
then $A_0=B^{\alpha/\beta}$ and $A_1=A$ have same order $\alpha$.  
The powers of $B$ can be defined as in Lemma \ref{power} using any cut 
$L_\theta=\{re^{i\theta},r\in(0,\infty)\}$ with $\theta\not=0(\pi)$
since $B$ is invertible self-adjoint with positive principal symbol,
moreover we have from Lemma \ref{power} that $\eta\in\Psi_\reg^{0,0}(M)=\Psi_\odd^{0,0}(M)$ 
using assumptions on $A,B$. The power $\eta^t$ and $\log \eta$ are constructed in the proof
of Corollary 2.1 in \cite{KV}: since $\eta$ is invertible with positive principal symbol and 
conjugate to the self-adjoint operator $B^{-\frac{\alpha}{2\beta}}AB^{-\frac{\alpha}{2\beta}}$  
its spectrum is then included in $[-R,-\eps]\cup [\eps,R]$ for some $0<\eps<R$,
one can thus take a cut $L_\theta$ for $\theta\not=0(\pi)$ to define the $\log$ and set 
\begin{equation}\label{etat}
\eta^t=-(2\pi i)^{-1}\int_{\Lambda_{\eps,R,\theta}}z^t(\eta-z)^{-1}dz
\end{equation}
where $\Lambda_{\eps, R,\theta}=\cup _{i=1}^4\Lambda^i$ 
\[\Lambda^1:\{re^{i\theta}, R\geq r\geq \eps\}, \quad \Lambda^2:=\{re^{i(\theta-2\pi)},\eps\leq r\leq R\},\quad
\Lambda^3:=\{\eps e^{i\varphi}, \theta\geq\varphi\geq \theta-2\pi\}\]
and $\Lambda^4$ is the circle of radius $R$ oriented opposite to the clockwise $(\theta-2\pi\leq \varphi\leq \theta)$.
The term $\log \eta$ is the derivative at $t=0$ of $\eta^t$ and this amounts to replace $z^t$ by $\log z$
in (\ref{etat}).
Using Corollary 2.1 in \cite{KV} the multiplicative anomaly is 
\begin{equation}\label{logfab}
\log F(A,B)=-\int_{0}^1 \textrm{WRes}_0\Big[(\log \eta)\Big(
\frac{\log (\eta^tB^{\frac{\alpha+\beta}{\beta}})}{\alpha+\beta} -
\frac{\log (A_t)}{\alpha}\Big)\Big]dt\end{equation}
if $F(A,B)=\det (AB)(\det A\det B)^{-1}$ and $\textrm{WRes}_0$ is the first Wodzicki residue
defined in (\ref{wod}).
From \cite[Prop. 4.2]{KV}, we know that $\eta^t$ and $\log \eta$ are in $\Psi_\odd^{0,0}(M)$
thus $A_t\in\Psi^{\alpha,0}_\reg(M)$ since $B^{\alpha/\beta}$ has regular parity by Lemma
\ref{power}. It remains to observe that the proofs of \cite[Prop. 4.3]{KV} and \cite[Cor. 4.3]{KV}
go through as well if we deal with non-integer $\psdo$'s with regular parity 
(as long as their order and sum of orders are non-zero), 
which is quite direct and essentially done using the same arguments than 
in Lemma \ref{power}; this proves that 
\[\frac{\log (\eta^tB^{\frac{\alpha+\beta}{\beta}})}{\alpha+\beta} -
\frac{\log (A_t)}{\alpha}\in\Psi^{0,0}_\odd(M)\] 
and after multiplying by $\eta\in\Psi_\odd^{0,0}(M)$ 
the Wodzicki residue in (\ref{logfab}) vanishes. The proof is achieved.
\qed\\

As a consequence we get the formulae
\begin{eqnarray*}
\det (P^{-\frac{\la}{2}+\frac{n}{4}}S(\la)P^{-\frac{\la}{2}+\frac{n}{4}+1})&=& \det (P^{-\frac{\la}{2}+\frac{n}{4}}S(\la)P^{-\frac{\la}{2}+\frac{n}{4}})\det P \\
&=&\det (P^{-\frac{\la}{2}+\frac{n}{4}})\det (S(\la))\det (P^{-\frac{\la}{2}+\frac{n}{4}+1}).
\end{eqnarray*}
Using again Lemma \ref{multdet}, we have 
\[\det (P^{-\frac{\la}{2}+\frac{n}{4}})\det (P^{-\frac{\la}{2}+\frac{n}{4}+1})(\det P)^{-1}=
\det (P^{-\la+\ndemi})\] 
if $\la\not=\{\ndemi,\ndemi+2\}$. We thus conclude, in view of definition (\ref{dettils}) that 
\begin{equation}\label{relationsstil}
\det \til{S}(\la)=\det(P^{-\la+\ndemi})\det S(\la)
\end{equation}
Remark that it is straightforward to see that 
\begin{equation}\label{detpa}
\det (P^\alpha)=e^{\alpha \textrm{TR}(\log P)}, \quad \alpha \in \rr^*
\end{equation}
using exactly same arguments than those we explained for $S(\la)$.\\ 
 
Now the interesting fact is that $\til{S}(n-\la)\til{S}(\la)=1$, and $\det(1)=\det(P)/\det(P)=1$
as defined by the method of \cite[Cor. 4.1]{KV}, thus 
\begin{eqnarray*}
1&=&\det(1)=\det(\til{S}(\la)\til{S}(n-\la))=\frac{\det(\til{S}(\la)\til{S}(n-\la)P^2)}{\det(P^2)}\\
&=&\frac{\det(\til{S}(\la)P)}{\det P}\frac{\det(\til{S}(n-\la)P)}{\det P}=\det\til{S}(\la)\det\til{S}(n-\la)
\end{eqnarray*} 
where we used (\ref{propdet}) and Lemma \ref{multdet} for $\det(P^2)=(\det P)^2$.
This equation can be rephrased, in view of (\ref{relationsstil}) and 
(\ref{detpa}), into 
\[\det S(\la)\det S(n-\la)=1.\]
Combining this identity with Proposition \ref{mainth} and the fact that $\xi$ is odd, we finally deduce  
\begin{equation}\label{plusmoins}
\det S\Big(\ndemi+iz\Big)=\pm e^{-2i\pi\xi(z)}\end{equation}
where $\pm$ means that one can not decide if it is plus or minus at the moment. 
As a matter of fact, it is proved in \cite[Lem. 4.16]{PP} or \cite[Lem. 4.3]{GZ}
that $S(\ndemi)=1-2P_X$
for some finite dimensional projector $P_X$ with rank equal to the multiplicity $m(n/2)$ of $n/2$ as 
resonance of $\Delta_g$. Using Proposition 6.4 of \cite{KV} we obtain 
\[\det S(\ndemi)=\det \til{S}(\ndemi)=\frac{\det (S(\ndemi)P)}{\det P}=\textrm{det}_{\textrm{Fr}}S(\ndemi)\]
where $\det_{\textrm{Fr}}$ is the Fredholm determinant defined for operators of 
the form ``Identity plus trace class''. But here it is trivial to see, from properties of $S(n/2)$ that 
$\det_{\textrm{Fr}}S(n/2)=(-1)^{m(\ndemi)}$. As a conclusion, the sign in (\ref{plusmoins}) is $(-1)^{m(\ndemi)}$ 
since $\xi(0)=0$. 
We bring together all these informations in the  
\begin{theo}\label{dettilfinal}
If $(X,g)$ is an even dimensional asymptotically hyperbolic manifold even modulo $O(x^\infty)$
and $x$ a geodesic boundary defining function. The 
generalized determinant of $S(\la)$ is
\[\det S\Big(\ndemi+iz\Big)=(-1)^{m(\ndemi)} e^{-2i\pi \xi(z)},\]
$\xi(z)$ being the extension, defined modulo $\zz$, of the generalized Krein function.
In particular if $n^2/4-k^2\notin\sigma_{\textrm{pp}}(\Delta_g)$, we have   
\[\det P_k=(-1)^{m(\ndemi)} e^{-2i\pi\xi(-ik)}\]
if $P_k$ is invertible whereas $\det P_k=0$ if $\ker P_k\not=0$.
\end{theo}

Notice that the determinant $\det P_k$ is a conformal invariant of the conformal infinity $(M,[h_0])$
of $(X,g)$ which depends only on the first $2k$ derivatives $(\pl_x^{j}x^2g|_{M})_{j\leq 2k}$. 
In the case of an even Poincar\'e-Einstein manifold, this is the determinant
of the $k$-th GJMS conformal Laplacian of $(M,[h_0])$.

It is clear that $\det P_k=0$ if $\ker P_k=\ker S(n/2+k)\not=0$
and $n^2/4-k^2\notin\sigma_{\textrm{pp}}(\Delta_g)$, 
indeed we can use Proposition \ref{mainth} to see that the divisor of $\det S(\la)$ 
is positive at $n/2+k$.\\
 
\subsection{Application to Selberg zeta function}
We conclude this discussion by an application to convex co-compact hyperbolic manifolds.
Let us first define the function
\begin{equation}\label{lt}
L(t)=\frac{\Gamma(\ndemi+it)\Gamma(\ndemi-it)}{\Gamma(it)\Gamma(-it)}=\Big(\frac{1}{4}+t^2\Big)\dots
\Big(\big(\ndemi-1\big)^2+t^2\Big)t\tanh(\pi t)
\end{equation} 
then we obtain by integrating $\pl_z\xi(z)$ and using Theorem \ref{dettilfinal}, Corollary
\ref{convexcc}  
\begin{theo}\label{detpkzeta}
Let $X=\Gamma\backslash\hh^{n+1}$ be a convex co-compact quotient of even dimension of $\hh^{n+1}$, let $S_X(\la),S_{\hh^{n+1}}(\la)$
be the respective scattering operator of $X$ and $\hh^{n+1}$. Then
\begin{equation}\label{detsn} 
\frac{Z(\ndemi-iz)}{Z(\ndemi+iz)}= \frac{\det S_X\Big(\ndemi+iz\Big)}{\Big(\det S_{\hh^{n+1}}\Big(\ndemi+iz\Big)\Big)^{\chi(\bar{X})}}=
\det S_X\Big(\ndemi+iz\Big)\exp\Big(
\frac{2i\pi(-1)^{\frac{n+1}{2}}}{\Gamma (n+1)}\chi(\bar{X})\int_{0}^z L(t)dt\Big)
\end{equation}
where $\chi(\bar{X})$ is the Euler characteristic of $\bar{X}$, $Z(s)$ the Selberg zeta function of 
the group $\Gamma$ and $L(t)$ is defined in \eqref{lt}. If $P_k$ is the GJMS $k$-th conformal Laplacian of its conformal infinity, 
if $n^2/4-k^2\notin\sigma_{\textrm{pp}}(\Delta_X)$ and if $P_k$ is invertible we have
\begin{equation}\label{detpkf}
\det P_k=\frac{Z(\ndemi-k)}{Z(\ndemi+k)}\exp\Big(\frac{2\pi(-1)^{\frac{n+3}{2}}}{\Gamma(n+1)}\chi(\bar{X})\int_{0}^kL(-it)dt\Big)\end{equation}
The integrals are understood as contour integrals avoiding the singularities, the final result remaining independent as proved before.
\end{theo}
\textsl{Proof}: we integrate $-2\pi i\pl_z\xi(z)$ using Corollary
\ref{convexcc} and the fact that $Z(\la)$ is holomorphic with no zeros on $\{\Re(\la)=\ndemi; \la\not=\ndemi\}$ (see \cite[Th. 6.2]{PP})
\[e^{-2i\pi (\xi(z)-\xi(\eps))}=\frac{Z(\ndemi-iz)Z(\ndemi+i\eps)}{Z(\ndemi+iz)Z(\ndemi-i\eps)}\exp\Big(-i\pi^{-\ndemi}
\frac{\Gamma(\ndemi)}{\Gamma (n)}\zerov(X)\int_{\eps}^z L(t)dt\Big) \]
for $\eps>0$ small and $z>\eps$.
Now we let $\eps\to 0$ and use the fact that $Z(\la)$ has a zero of order $m(\ndemi)$ at $\la=\ndemi$ (see \cite[Th. 6.2]{PP})
to see that 
\[\lim_{\eps\to 0}\frac{Z(\ndemi+i\eps)}{Z(\ndemi-i\eps)}=(-1)^{m(\ndemi)}\]
and since $L(t)$ is regular at $t=0$ and $\xi(0)=0$, 
Theorem (\ref{dettilfinal}) gives the result by meromorphic continuation in $z$. Note that 
we used the identity $\zerov(X)=(-1)^{\frac{n+1}{2}}\pi^{\ndemi+1}\chi(\bar{X})/\Gamma(\ndemi+1)$
and that the value $\det S_{\hh^{n+1}}(\la)$ is just obtained by taking the trivial group 
$\Gamma=\{\rm Id\}$, and $Z(\la)=1$ in that case in Corollary \ref{convexcc}.
\qed\\

It is interesting to compare such a result with that of Sarnak \cite{Sa} for instance, where he 
proved for Riemann surfaces the identity between some determinant of Laplacian with the Selberg 
zeta function, see also the recent paper of Borthwick-Judge-Perry \cite[Th. 5.1]{BJP2} where 
in this case this is the determinant of the Laplacian on the interior (the non-compact 
hyperbolic manifold) which is related to $Z(\la)$. An approach as in \cite{BJP}
could also be used to define a generalized determinant $d(z)=\det(\Delta_g-n^2/4-z)$, our results 
would be interpreted as an identity relating $d(z+i0)/d(z-i0)$ to $\det S(n/2+iz^\demi)$ for $z\in(0,\infty)$.\\

\subsection{The Weyl asymptotic}
Finally we remark that a Weyl type asymptotic holds for $\xi$ for hyperbolic 
convex co-compact quotient.
\begin{prop}\label{weyl}
If $X=\Gamma\backslash\hh^{n+1}$ is a convex co-compact quotient with $n$ odd and 
$\delta$ is the dimension of the limit set of $\Gamma$, then as $t\to\infty$
\[\xi(t)=
\left\{\begin{array}{ll}
\frac{(4\pi)^{-\frac{n+1}{2}}}{\Gamma(\frac{n+3}{2})}\zerov (X)\Big(t^{n+1} +\sum_{i=1}^{\frac{n-1}{2}}C_it^{2i}\Big)+O(t), & {\rm if }\delta<\ndemi \\
\frac{(4\pi)^{-\frac{n+1}{2}}}{\Gamma(\frac{n+3}{2})}\zerov (X)t^{n+1} +O(t^n), & \rm{ otherwise}
\end{array}\right.\]
where $C_i$ is the $t^{2i}$ coefficient of the polynomial
\[\int_{0}^tu\prod_{j=1}^{\frac{n-1}{2}}(\ndemi-j+u^2)du.\]
\end{prop}
\textsl{Proof}: For the case $\delta<\ndemi$ we can apply Corollary \ref{convexcc} with the estimate on zeta function
\[\left|\frac{Z'(\la)}{Z(\la)}\right|+\left|\frac{Z'(n-\la)}{Z(n-\la)}\right|\leq C , \quad \Re(\la)=\ndemi\]
clearly obtained from its definition in (\ref{zeta}) and the absolute convergence of this sum on
the line $\Re(\la)=n/2$ if $\delta<n/2$.  
For the general case, it suffices to apply the same proof as that of 
Melrose \cite{ME0}, see also \cite{GZ} for Riemann surfaces. 
The only ingredients needed are a representation of $\pl_z\xi$
\[\pl_z\xi(z)= (2\pi)^{-1}\sum_{\la \in \mc{D}}\Big(\frac{2\Im(\la)}{(z-\Re(\la))^2+\Im(\la)^2}
+P_{\la}(z)\Big)+P(z), \quad z \in\rr\] 
for $\mc{D}\subset\cc\setminus \rr$ a discrete set such that $\{\la\in\mc{D}; |\la|<R\}=O(R^{n+1})$ (counted with multiplicites) with $P_\la(z),P(z)$ polynomials in
$z$ of degree $\leq n$, and the singularity of its Fourier transform (i.e. 
the $0$-Trace of the wave kernel) at $t=0$ with first coefficient 
given by $\zerov(X)$ times the appropriate constant. The representation of $\pl_z\xi$ is a trivial consequence
of a Hadamard factorization of $\det S(\la)$ as a quotient of two entire functions 
of order $n+1$ with symmetric respective zeros given by the sets $\mc{D}$ and 
$\bar{\mc{D}}=\{\la; \bar{\la}\in\mc{D}\}$. This factorization is clearly deduced from two facts:
first the analysis of the zeros and poles of $S(\la)$ in Proposition \ref{mainth}, with 
$\mc{D}$ being the set of resonances combined with the $\{n/2+k\in n/2+\nn; \ker S(n/2+k)\not=0\}$ counted
with multiplicities, the symmetry of the poles of $S(\la)$ with respect to $\la\to \bar{\la}$
comes classically from $S(n-\bar{\la})^*=S(\la)^{-1}=S(n-\la)$ (see \cite{GRZ}); 
secondly the fact that $\det S(\la)$ is meromorphic of order $n+1$, which by (\ref{detsn})  
is a consequence of the Hadamard factorisation of $Z(\la)$ proved in \cite[Th. 1.1.]{PP}
(note that the growth of the functions are obtained by thermodynamic formalism 
of \cite{Fr,B}) and a Hadamard factorisation of 
\[M(z):=\exp\Big(-i\pi^{-\ndemi}\frac{\Gamma(\ndemi)}{\Gamma (n)}\zerov(X)\int_{0}^z L(t)dt\Big)\]
as a quotient of two entire functions of order $n+1$, that we need to check.  
To prove it, it suffices to remark that its zeros and poles are at $z=\pm i(n/2+k)$ ($k\in\nn_0$)
with multiplicity $|\chi(\bar{X})|h_n(k)<C|k|^{n}$ where $C$ is a constant and $h_n(k)$
is the dimension of the space of spherical harmonics of degree $k$ on $S^{n+1}$ (see \cite[Rem. 6.10]{PP}), 
we also have a trivial bound $|M(z)|<e^{C'(|z|+1)^{n+1}}$ in $\textrm{dist}(z, \pm i(n/2+\nn))>1/4$ 
for some constant $C'$ by estimating $L(t)$ in (\ref{lt}), thus by multiplying by a Hadamard product
of order $n+1$ that has zeros at poles of $M(z)$ and using the maximum principle, 
we have our proof. 

The singularity of the wave $0$-Trace is studied by Joshi-Sa Barreto \cite{JSB2}, 
thus we are done by applying Melrose's proof \cite{ME0}.
\qed\\

\section{Appendix - Proof of Proposition \ref{kvtrace}}

The part with $\supp \chi_i\cap\supp\chi_j=\emptyset$ is clear, the other part can be worked 
out in each $U_{ij}\x U_{ij}$ using polar coordinates (or blow-up coordinates), and it clearly suffices
to assume that the kernel of $F,L$ are supported in $U_{ij}\x U_{ij}$, by partition of unity arguments. 
The volume density on $U_{ij}$ is trivialized by $|\det h_0|^\demi dy'$, 
let us include $\chi_i\otimes(\chi_j|\det h_0|^\demi)$ into $F$ to simplify notations: 
$F(y,y')$ in $U_{ij}\x U_{ij}$ will now mean $F(y,y')\chi_i(y)\chi_{j}(y')|\det h_0(y')|^{\demi}$.
The odd parity property of $F$ still holds in view of
properties of odd classes discussed before.
We denote by $U_{ij}\x_0 U_{ij}=\beta^{-1}(U_{ij}\x U_{ij})$ where 
\[\beta: U_{ij}\x [0,\infty)\x S^{n-1}\to U_{ij}\x \rr^{n}, \quad \beta(y,r,\omega)=
(y,y+r\omega).\]
By assumption on $F,L$, the lifted kernel of $F,L$ under $\beta$ in $U_{ij}\x U_{ij}$ 
can be decomposed as (using that $L$ is symmetric)
\[\beta^*F(y,r,w)=F(y,y')=F_{\sing}(y,r,w)+F_{\reg}(y,r,w),\]
\[\beta^*L(y,r,w)=L(y',y)=L_{\sing}(y,r,w)+L_{\reg}(y,r,w),\]
\[F_{\sing}(y,r,w)=\sum_{i=0}^n(F_i(y,w)+F_{i,\log}(y,w)\log r)r^{-2\la+i},\quad
L_{\sing}(y,r,w)=\sum_{i=0}^nL_i(y,w)r^{-2n+2\la+i}\]
where $r=|y-y'|, w=(y'-y)/r$,  
and finally $F_i,F_{i,\log},L_i$ smooth satisfying
\begin{equation}\label{transm}
F_i(y,-w)=(-1)^iF_i(y,w), \quad F_{i,\log}(y,-w)=(-1)^iF_{i,\log}(y,w), \quad
L_i(y,-w)=(-1)^i L_i(y,w)\end{equation}
whereas $F_{\reg}\in\Psi^{2\la-2n-1,1}(U_{ij})$, $L_{\reg}\in\Psi^{-2\la-1,0}(U_{ij})$. Note that 
the distributions $\beta_*L_\sing$ and $\beta_*L_\reg$ can be extended in $M\x M$ 
so that $\beta_*L_\reg+\beta_*L_\sing$ is the kernel of $L$ (the extension is not 
relevant for what follows but simplifies statements). 
Since $F$ has compact support in $U_{ij}\x U_{ij}$ it is possible to take constants $B>A>0$ such that 
\[\forall A'\in[A,B], \supp(\beta^*F)\subset \Big((U_{ij}\x [0,A']\x S^{n-1})\cap \supp(\beta^*F)\Big)\subset U_{ij}\x_0 U_{ij}\]
and define $\psi\in C_0^\infty([0,B))$ equal to $1$ on $[0,A]$ then 
$\psi(r)\beta^*F=\beta^*F$.

Identifying operators and Schwartz kernels, the $3$ composed operators 
$\beta_*(\psi(r) F_{\sing})\beta_*(L_{\reg})$, 
$\beta_*(\psi(r) F_{\reg})\beta_*(L_{\sing})$, 
$\beta_*(\psi(r) F_{\reg})\beta_*(L_{\reg})$ are operators 
in $\Psi^{-n-1,1}(M)$ thus of trace class, their KV-Trace 
is the usual trace, that is the integral on the diagonal or equivalently, using polar cordinates,
\[\TR(\beta_*(\psi(r) F_{\sing})\beta_*( L_{\reg}))=\int_{U_{ij}}\int_0^\infty\int_{S^{n-1}}
\psi(r) F_{\sing}(y,r,\omega)L_{\reg}(y,r,\omega)d\omega dr\textrm{d}_{h_0}(y)\]
and the obvious similar formula for both other operators. 

Let us now deal with the singular 
term $\beta_*(\psi F_{\sing})\beta_*(L_{\sing})$, according to corollary \ref{corol}, we have to study the kernel 
\[M(y,z)=\int \psi(|y-y'|)F_{\sing}\Big(y,|y'-y|,\frac{y'-y}{|y'-y|}\Big)
L_{\sing}\Big(z,|y'-z|,\frac{y'-z}{|y'-z|}\Big)dy'\]
near $y=z$, or to consider $M(y,y+u)$ for $u\in\rr^n$ small,
make an expansion in homogeneous functions of $u$ near $u=0$ and get rid of the divergent terms.

We begin with the following term for $i+j<n$ and small $u$ 
\[M_{i,j,\log}(y,u):=\int \frac{\psi(|y-y'|)\log(|y'-y|)}{|y'-y|^{2\la-i}|y'-y-u|^{2n-2\la-j}}
F_{i,\log}\Big(y,\frac{y'-y}{|y'-y|}\Big)
L_{j}\Big(y+u,\frac{y'-y-u}{|y'-y-u|}\Big)dy'\]
\[\quad\quad\quad\quad\quad\quad=\int \psi(|y'|)\log(|y'|)|y'|^{-2\la+i}|y'+u|^{-2n+2\la+j}
F_{i,\log}\Big(y,\frac{y'}{|y'|}\Big)
L_{j}\Big(y+u,\frac{y'+u}{|y'+u|}\Big)dy'.\]
We split the integral in two $M_{i,j,\log}^{(1)},M_{i,j,\log}^{(2)}$ 
by inserting cut-off functions $(1-\chi)(|y'|/|u|)$ and $\chi(|y'|/|u|)$  
where $\chi\in C_{0}^\infty([0,1))$ equals to $1$ in $[0,1/2]$. Since $\chi(|y'|/|u|)\psi(|y'|)=\chi(|y'/|u|)$
for small $|u|$, we obtain through a change of variable $y'\to y'|u|$ 
that the first integral is, after having set $u=|u|\theta$,
\[M_{i,j,\log}^{(1)}(y,u)=|u|^{-n+i+j}\int \chi(|y'|)\log(|y'|)\frac{|y'|^{-2\la+i}}{|y'+\theta|^{2n-2\la-j}}
F_{i,\log}\Big(y,\frac{y'}{|y'|}\Big)
L_{j}\Big(y+u,\frac{y'+\theta}{|y'+\theta|}\Big)dy'\]
\[\quad\quad\quad\quad\quad\quad+|u|^{-n+i+j}\log(|u|)\int \chi(|y'|)\frac{|y'|^{-2\la+i}}{|y'+\theta|^{2n-2\la-j}}
F_{i,\log}\Big(y,\frac{y'}{|y'|}\Big)
L_{j}\Big(y+u,\frac{y'+\theta}{|y'+\theta|}\Big)dy'.\]
Since for $|u|>0$ these are the values of the distributions 
$|y'|^{-2\la+i}$ and $|y'|^{-2\la+i}\log|y'|$ again smooth compactly supported functions 
in $y'$, these well-defined integrals for $\Re(\la)>\ndemi$ can be classically 
defined for $\la\notin\ndemi-\demi\nn$ by holomorphic extension in $s$ at $s=0$ after 
multiplying it by $|y'|^{s}$ for $\Re(s)\gg 0$
using Taylor expansion of the compactly supported function at $y'=0$. 
By differentiating under the integral one clearly gets an expansion as $|u|\to 0$ of the form 
\[M_{i,j,\log}^{(1)}(y,u)\sim |u|^{-n+i+j}\sum_{k=0}^\infty|u|^k(\alpha_k(y,\theta)+\log|u| \beta_k(y,\theta))\]
with $\alpha_k,\beta_k$ smooth. This is trivial if $\Re(\la)>\ndemi$ since the integral converges but not much more complicated if $\la\notin \ndemi-\demi\nn$ using the way the analytic extension in $s$ at $s=0$ is constructed. 
If $\Re(\la)>\ndemi$, each $\alpha_k,\beta_k$ is expressed as an integral and it is 
easy to see by change of variable $y'\to-y'$ in that integral and regular parity of $F,L$ that $\alpha_k(y,-\theta)=(-1)^{i+j+k}\alpha_k(y,\theta)$ and the same for $\beta_k$ but
$|u|^{-n+i+j+k}\alpha_k$ and $|u|^{-n+i+j+k}\beta_k$ are homogeneous of degree $-n+i+j+k$ which implies 
that $M_{i,j,\log}^{(1)}$ does not contribute to the KV-Trace using Corollary \ref{corol}; 
the same holds for any $\la\not\in\ndemi-\demi\nn$ by holomorphic continuation arguments in $s$ at $s=0$ 
(we let it as an exercise).
We consider the second term 
\[M_{i,j,\log}^{(2)}(y,u)=|u|^{-n+i+j}\int (1-\chi(|y'|))\frac{\log(|y'||u|)|y'|^{-2\la+i}}{|y'+\theta|^{2n-2\la-j}}
F_{i,\log}\Big(y,\frac{y'}{|y'|}\Big)
L_{j}\Big(y+u,\frac{y'+\theta}{|y'+\theta|}\Big)dy'\]
\[\quad\quad\quad\quad+\int (\psi(|y'|)-1)\frac{|y'|^{-2\la+i}\log|y'|}{|y'+u|^{2n-2\la-j}}
F_{i,\log}\Big(y,\frac{y'}{|y'|}\Big)
L_{j}\Big(y+u,\frac{y'+u}{|y'+u|}\Big)dy'\]
where we used $\psi(|y'|)\chi(|y'|/|u|)=\chi(|y'|/|u|)$ and a 
change of variable $y'\to |u|y'$ as before in the first integral. 
The first integral can be dealt with like $M_{i,j,\log}^{(1)}$ since now this is the value 
of the distribution $|y'+\theta|^{2n-2\la-j}$ again a smooth compactly supported function plus 
the integral of an $L^1(\rr^n,dy')$ function ($|y'|^{-2\la+i}|y'+\theta|^{2n-2\la-j}\log|y'|$ is $L^1$ 
on $\{y'\in\rr^n; |y'|>2\}$). Thus the first integral does not contribute to the KV-Trace.
Now the second integral is a smooth function of $u$ since $i+j<n$ and $1-\psi(|y'|)=0$ near $y'=0$ 
thus its KV-Trace is the integral of its value on the diagonal $\{u=0\}$, that is
\[\TR(M_{i,j,\log})=\int (\psi(|y'|)-1)|y'|^{-2n+i+j}\log|y'|
F_{i,\log}\Big(y,\frac{y'}{|y'|}\Big)L_{j}\Big(y,\frac{y'}{|y'|}\Big)dy'dy\]
where we used an abuse of notation by writing $\TR(M_{i,j,\log})$ to mean its contribution to the KV-Trace
of $FL$; equivalently
\[\TR(M_{i,j,\log})=\int_{U_{ij}}\int_0^\infty\int_{S^{n-1}} (\psi(r)-1)r^{-2n+i+j}\log r
F_{i,\log}(y,\omega)L_{j}(y,\omega)d\omega dr\textrm{d}_{h_0}(y).\]

The part $M_{i,j,\log}$ with $i+j>n$ gives a trace class operator, thus its KV-Trace is 
the trace, the integral of the kernel on the diagonal.

It remains to deal with the critical case where $i+j=n$. One can decompose as before with the 
function $\chi(|y'|/|u|)$ and same arguments show that the only term that could 
contribute is 
\[M_{i,j,\log}^{(2)}(y,u)=\int_{S^{n-1}}\int_{0}^{|u|^{-1}} (1-\chi(r))\frac{\log(r|u|)r^{-2\la+i+n-1}}{|r\omega+\theta|^{2n-2\la-j}}
F_{i,\log}(y,\omega)
L_{j}\Big(y+u,\frac{r\omega+\theta}{|r\omega+\theta|}\Big)drd\omega\]
\[\quad\quad\quad\quad+\int_{|y'|<1} (\psi(|y'|)-1)\frac{|y'|^{-2\la+i}\log|y'|}{|y'+u|^{2n-2\la-j}}
F_{i,\log}\Big(y,\frac{y'}{|y'|}\Big)
L_{j}\Big(y+u,\frac{y'+u}{|y'+u|}\Big)dy'\]
here $r=|y'|$ and $\omega=y'/r$. The second integral is smooth in $u$, and its value at $u=0$ is
$0$ by using a change of variable $y'\to-y'$ and the oddness properties (\ref{transm}), thus it does not contribute to the KV-Trace. 
A Taylor expansion of $L_{j}(y+u,w)$ at $u=0$ induces an 
expansion of the first integral as $|u|\to 0$ of the form
\[ \sum_k |u|^k\int_{S^{n-1}}\int_0^{|u|^{-1}}\gamma_k(y,|u|,\theta,r,\omega)drd\omega\]
for some $\gamma_k$ where each integral is easily seen to be a $O(\log^2|u|)$ and 
\[\gamma_0(y,|u|,\theta,r,\omega)=(1-\chi(|r|))\frac{\log(r|u|)r^{-2\la+i+n-1}}{|r\omega+\theta|^{2n-2\la-j}}
F_{i,\log}(y,\omega)
L_j\Big(y,\frac{r\omega+\theta}{|r\omega+\theta|}\Big).\]
Remark by using change of variable $\omega\to-\omega$ and (\ref{transm})
that the integral of $\gamma_0 d\omega$ on the sphere (in variable $\omega$) is an odd function 
of $\theta$.
Making an expansion of the form $\sum_l\log(r|u|)r^{-1-l}\mu_l(y,\theta,\omega)$ 
as $r\to\infty$ of each $\gamma_k$  for some $\mu_l$ with $l\in\nn_0$, 
one gets an asymptotic expansion of $M_{i,j,\log}^{(2)}$ as $|u|\to 0$ of the form
\[M_{i,j,\log}^{(2)}(y,u)\sim c(y,u/|u|)+\log|u|\sum_{k=0}^\infty |u|^k\Big(\log|u|\alpha_k(y,u/|u|)+\beta_k(y,u/|u|)\Big)\]
for some $c,\alpha_k,\beta_k$ smooth and actually a small calculation gives 
\[-2\alpha_0(y,\theta)=\beta_0(y,\theta)=\int_{S^{n-1}}F_{i,\log}(y,\omega)
L_j(y,\omega)d\omega\]
which vanishes since the integrand is an odd function of $\omega$ by (\ref{transm}).
This implies that the function $M_{i,j,\log}^{(2)}(y,|u|\theta)$ has a limit as $|u|\to 0$ 
which is $c(y,\theta)$ but also 
\[c(y,\theta)=\lim_{|u|\to 0}\int_0^{|u|^{-1}}\int_{S^{n-1}}\gamma_0(y,|u|,\theta,\omega)drd\omega.\] 
It remains to observe that this is an odd function of $\theta$ since it is for any $u>0$.
Then the KV-Trace of $M_{i,j,\log}^{(2)}(y,|u|\theta)$ vanishes and the same arguments work as well for 
the part without $\log$ terms in the integral defining $M(y,y+u)$ (they are actually simpler).\\ 

So far we have proved that (recall that we have included $\chi_i,\chi_j$ in $F$) 
\[\TR(FL)=\int\Big(\psi F_\sing L_\reg+\psi F_\reg L_\sing+\psi F_\reg L_\reg
+\psi[F_\sing L_\sing]_\reg\Big)\quad\quad\]
\[\quad\quad\quad\quad\quad\quad\quad\quad\quad\quad+(\psi(r)-1)[\beta^*F
\beta^*L]_{\sing} dw_{S^{n-1}}dr\textrm{d}_{h_0}(y).\]
To conclude, we observe that 
\[ [\beta^*F\beta^*L]_{\sing}=[F_\sing L_\sing]_{\sing}\]
\[[\beta^*F\beta^*L]_\reg=(\psi-1)F_\sing L_\sing+\psi F_\sing L_\reg+\psi F_\reg L_\sing+\psi F_\reg L_\reg
+[F_\sing L_\sing]_\reg\]
and since the result holds for arbitrary $\psi$, 
this is also true for the limit case $\psi=\indic_{[0,A]}$. The desired formula holds for the chosen $A$. 
Now a straightforward computation leads to
\begin{equation}\label{TRFL}
\TR(FL)=\int \indic_{[0,A']}(r)[\beta^*F\beta^*L]_\reg-\indic_{[A',\infty)}(r)[\beta^*F\beta^*L]_\sing drd\omega_{S^{n-1}}\textrm{d}_{h_0}(y)
\end{equation}
\[\quad \quad\quad-\int \indic_{[A,A']}(r)\Big([\beta^*F\beta^*L]_{\reg}+[\beta^*F\beta^*L]_{\sing}\Big) drd\omega_{S^{n-1}}\textrm{d}_{h_0}(y)\] 
if $A'>A$. Since $[\beta^*F\beta^*L]_\reg+[\beta^*F\beta^*L]_\sing=\beta^*F\beta^*L$ 
up to the $n$-th homogeneous 
(and log-homogenous) term and since that term has vanishing integral (it is odd in the sphere variable) 
on $U_{ij}\x [A,A']\x S^{n-1}$, the second integral in $(\ref{TRFL})$ is the integral of $\beta^*F\beta^*L$ whose support
is included in $r\in[0,A]$, thus it vanishes and we have proved the desired formula 
for any $A'>A$. 
\qed\\

\textbf{Acknowledgements}: First I would like to thank Pierre Albin for interesting discussions 
on renormalization of the $0$-Trace and Peter Perry for sharing his approach on the subject. 
Tom Branson, Gilles Carron, Frederic Naud, Paul Loya and Martin Olbrich 
have also been helpful by answering my questions
and Laurent Guillop\'e, William Ugalde for their comments. 
This work has begun at NCTS center of Hsinchu (Taiwan), I am grateful to 
the organizers -and particularly Alice Chang- of the ``conformal invariants'' meeting for their invitation.
It has then be continued in the mathematics department of Nice and completed
at the mathematics department of ANU, Canberra (Australia), during a research associateship funded by the Australian  
Research Council. It is partially supported by french 
ANR grant JC05-52556 and NSF grant DMS0500788.

\end{document}